\documentclass[12pt]{amsart}

\usepackage{amsfonts, amssymb, amsmath, latexsym,  mathtools, amscd}
\usepackage{hyperref,graphicx}
\usepackage[dvipsnames]{xcolor}
\usepackage{url}
\usepackage{enumitem}
\usepackage{wrapfig}
\usepackage{ulem}  % rule out text

\hypersetup{
    colorlinks=true,
    linkcolor=black,
    filecolor=magenta,      
    urlcolor=blue,
    }
\newtheorem{theorem}{Theorem}[section]

\newtheorem{lemma}[theorem]{Lemma}

\newtheorem{proposition}[theorem]{Proposition}

\usepackage{algorithm}
\usepackage{algpseudocode}

\makeatletter
\let\c@algorithm\c@theorem

\makeatother

\theoremstyle{remark}

\newtheorem{definition}[theorem]{Definition}
\newenvironment{Definition}
{\pushQED{\qed}\definition}
{\popQED\enddefinition}

\newtheorem{remark}[theorem]{Remark}
\newenvironment{Remark}
{\pushQED{\qed}\remark}
{\popQED\endremark}

\newtheorem{example}[theorem]{Example}
\newenvironment{Example}
{\pushQED{\qed}\example}
{\popQED\endexample}

\newcommand{\CC}{{\mathbb C}}
\newcommand{\NN}{{\mathbb N}}
\newcommand{\PP}{{\mathbb P}}
\newcommand{\QQ}{{\mathbb Q}}
\newcommand{\RR}{{\mathbb R}}
\newcommand{\TT}{{\mathbb T}}
\newcommand{\ZZ}{{\mathbb Z}}

\newcommand{\calB}{{\mathcal B}}
\newcommand{\calE}{{\mathcal E}}
\newcommand{\calF}{{\mathcal F}}
\newcommand{\calN}{{\mathcal N}}

\newcommand{\calV}{{\mathcal V}}

\newcommand{\BD}{{\mathrm{BD}}}
\newcommand{\rec}{{{\rm rec}}}
\newcommand{\disc}{{\mathrm{disc}}}
\newcommand{\vol}{{\mathrm{vol}}}
\newcommand{\CP}{{\mathrm{CP}}}
\newcommand{\Tr}{{\mathrm{Tr}}}
\newcommand{\Hess}{{\mathrm{Hess}}}
\def\cprime{$'$}

\newcommand{\defcolor}[1]{{\color{blue}#1}}
\newcommand{\demph}[1]{\defcolor{{\sl #1}}}

\title[Bloch Discriminants]{Bloch Discriminants of Discrete Periodic Operators}
\author{Margaret H.~Regan}
\address{Margaret H.~Regan, Department of Mathematics, Grand Valley State University, Allendale, MI 49424, USA}
\email{reganmar@gvsu.edu}
\urladdr{https://www.margarethregan.com}
\author{Frank Sottile}
\address{Frank Sottile, Department of Mathematics,
         Texas A\&M University, College Station, Texas 77843,  USA}
\email{sottile@tamu.edu}
\urladdr{https://franksottile.github.io/}
\author{Simon Telen}
\address{Simon Telen, Max Planck Institute for Mathematics in the Sciences,
      Inselstra{\ss}e 22, 04103 Leipzig, Germany}
\email{simon.telen@mis.mpg.de}
\urladdr{https://sites.google.com/view/simontelen/home-page}
\thanks{Research of Sottile supported in part by NSF grant DMS-2201005 and Simons TSM-00014009.}
\thanks{This research was started at the 2023 ICERM Workshop on Algebraic Geometry in Spectral Theory.}
\begin{document}
\begin{abstract}
  The Bloch variety of an operator on a periodic graph encodes the relation between its spectrum and the group action on the graph,
  with critical points corresponding to edges of the spectrum.
  We adapt methods from computational algebraic geometry to describe how the qualitative behavior of critical points of
  Bloch varieties depends upon parameters for three $\ZZ^2$-periodic graphs. 
  The set of parameters where the qualitative behavior changes is the Bloch discriminant.  
  As the critical points are real-algebraic with respect to a non-standard real structure on the complex torus, this is
  also a first step towards developing computational tools for varieties with a non-standard real structure.
\end{abstract}
%%%%%%%%%%%%%%%%%%%%%%%%%%%%%%%%%%%%%%%%%%%%%%%%%%%%%%%%%%%%%%%%%%%%%%%%%%%%%%%%%%%%%%%%%%%%%%%%%%%%

\maketitle

\section*{Introduction}

Discrete operators on periodic graphs are an approximation to physical models of crystals, e.g., the tight binding
model in solid state physics~\cite{AshcroftMermin,Kittel}.
The spectrum of the operator is the set of energies at which electrons may propagate through the medium.
After discretizing, the spectra of such operators may be studied from the perspective of (real) algebraic geometry via
objects such as the Bloch and Fermi varieties associated to it.  
This point of view was taken in~\cite{DKS}, which considered how the behavior of Bloch varieties at extrema depends
upon free parameters, observing in Theorem 12 {\it loc.\ cit.}\ that the Bloch variety has degenerate critical points either
at almost all parameter values or at almost no parameter values.
This algebraic perspective is developed in~\cite{AAPGO}.
We extend that from the vantage of computational algebraic geometry.
We aim to describe explicitly how the features of a Bloch variety depend upon parameters.
We develop some general theory and work out case studies for several graphs, which 
include the hexagonal lattice and a graph whose Bloch variety has a curve of degenerate critical points, similar
to the example of Filonov and Kachkovskiy~\cite{FK18}.

Let $\Gamma$ be a graph equipped with a free action of $\ZZ^d$ having finitely many orbits on its edges $\calE$ and vertices $\calV$.
We label $\Gamma$ with parameters $p=(e,V)$ where $e\colon\calE\to \RR^\times$ and $V\colon\calV\to\RR$ are constant on orbits.
These data define a
Schr\"odinger operator $H_p$ on a Hilbert space $\ell_2(\Gamma)$.
Applying the Floquet transform (a version of the Fourier transform), $H_p$ is represented by a matrix $H_p(z)$ of Laurent
polynomials in the Floquet parameters $z\in\TT^d$.
Here, $\TT\subset\CC^\times$ is the group of unit complex numbers.
For $z\in\TT^d$, the Floquet matrix $H_p(z)$ is Hermitian, and so all of its eigenvalues are real.
The Bloch variety $\calB_p$ is the hypersurface in $\TT^d\times\RR$ defined by the vanishing of the dispersion polynomial
$D_p(z,\lambda) = \det(\lambda\cdot{\rm id}-H_p(z))$.
Its projection to the $\lambda$-axis is the spectrum of $H_p$.
Fibers of this projection are Fermi varieties and critical points correspond to singular Fermi varieties.
We recall this in Section \ref{S:BlochVarieties}.

As $H_p(z^{-1}) = H_p(z)^\top$, the (complexified) Bloch variety is stable under the involution $z \mapsto z^{-1}$
and polynomials $F(z)$ defining Fermi varieties are reciprocal in that $F(z)=F(z^{-1})$.
Singular reciprocal hypersurfaces are parameterized by \demph{reciprocal $A$-discriminants} in analogy with classical
$A$-discriminants~\cite{GKZ}.
We carry out the first steps of their study in Section \ref{S:ReciprocalDiscriminants}.
A general theory of discriminants of reciprocal polynomials is developed in~\cite{RDisc}.

The \demph{Bloch discriminant} is the subset of the parameter space where a qualitative change of the critical points occurs.
These are critical points of the coordinate function $\lambda$ on the Bloch variety and we are concerned with their number,
location, multiplicity, etc., as the parameters $p=(e,V)$ vary.
The analysis of these critical points is of significant current interest~\cite{BCCM_LocalGlobal,DKS,FRS25,FS24}. 
The critical values are roots of the reciprocal $A$-discriminant evaluated at the coefficients of $z$ in $D_p(z,\lambda)$.
The complement of the Bloch discriminant consists of connected components called chambers with the property that
all operators whose parameters lie in a single chamber have qualitatively similar arrangements of critical points.
We make this precise in Definition~\ref{def:Bloch_Discriminant}. 
Describing the Bloch discriminant and the arrangement of critical points for each chamber is an instance of the
general problem of determining the geography of parameter space~\cite{BERNAL,Harrington_Decomp}.

Besides illuminating the behavior of the family of operators, this question is challenging for Bloch varieties, as the
features we are studying (critical points) are real points of $(\CC^\times)^d\times\CC$ for a non-standard real structure:
That is, the usual complex conjugation $(z,\lambda)\mapsto (\overline{z},\overline{\lambda})$ is replaced by the
antiholomorphic involution $(z,\lambda)\mapsto (1/\overline{z},\overline{\lambda})$.
The fixed (real) points of this new structure are $\TT^d\times\RR$, in contrast to the classical real points $(\RR^\times)^d\times\RR$.
A goal of this paper is to initiate the development of computational and analytical tools for real algebraic geometry of
algebraic varieties with a non-standard real structure.

This paper is organized as follows.
Section~\ref{S:BlochVarieties} defines Bloch varieties and their corresponding Bloch discriminants.
The $A$-discriminant of a reciprocal polynomial and its connection to the Bloch discriminant is discussed in
Section~\ref{S:ReciprocalDiscriminants}.
Section~\ref{Sec:bipartite} further investigates the connection of certain resultants and Hessians to the classification of the critical
points of the dispersion polynomial.
The paper closes with a discussion of these details for several bipartite graphs, including the hexagonal lattice, in
Sections~\ref{S:Hexagonal}, \ref{S:four}, and \ref{S:five}.

%%%%%%%%%%%%%%%%%%%%%%%%%%%%%%%%%%%%%%%%%%%%%%%%%%%%%%%%%%%%%%%%%%%%%%%%%%%%%%%%%%%%%%%%%%%%%%%%%%%%
\section{Bloch varieties} \label{S:BlochVarieties}

\begin{figure}[htb]
  \centering
   \begin{picture}(122,92)
     \put(0,0){\includegraphics{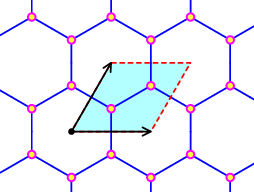}}
     \put(55,50){\small $W$}
   \end{picture}
   \qquad
   \begin{picture}(105,92)(0,-1)
     \put(0,0){\includegraphics{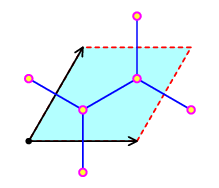}}
     \put(31,31){\small$u$}      \put(95,31){\small$u$}     \put(64,89){\small$u$}
     \put(69,53){\small$v$}      \put( 5,53){\small$v$}     \put(37,-2){\small$v$}
     \put(54,15){\small$x$}
     \put(29,62){\small$y$}
     \put(52,39){\small{\color{blue}$a$}}
     \put(32,14){\small{\color{blue}$c$}}   \put(68,73){\small{\color{blue}$c$}}
     \put(14,40){\small{\color{blue}$b$}}   \put(86,44){\small{\color{blue}$b$}}
   \end{picture}
   \caption{Hexagonal lattice with fundamental domain and a labeling.}
  \label{fig:HexGraph}
\end{figure}

We recall some notation and definitions from \cite{FS24}, which
may also be found in~\cite[Ch.\ 4]{BerkKuch}.
A $\ZZ^d$-periodic graph is a simple graph $\Gamma$ equipped with a free action of $\ZZ^d$ that has finitely many orbits on its vertices 
$\calV$  and edges $\calE$.
We  write the action of $\ZZ^d$ on $\calV$ additively: the image of $u\in\calV$ under
$\alpha\in\ZZ^d$ is $u{+}\alpha\in\calV$.
We choose a fundamental domain $\defcolor{W}\subset\calV$, which is a set of representatives of orbits of  $\ZZ^d$ on $\calV$.
Similarly, let \defcolor{$E$} be the set of $\ZZ^d$-orbits in $\calE$.
A labeling of $\Gamma$ is a pair $p=(e,V)$ of $\ZZ^d$-periodic functions $e\colon\calE \rightarrow \RR^\times$ and 
$V\colon \calV \rightarrow\RR$.
(An edge with label zero may be deleted to give a different graph.)
Equivalently, $p=(e,V)$ is a point in  $(\RR^\times)^E\times\RR^W$.
Figure~\ref{fig:HexGraph} shows (part of) the $\ZZ^2$-periodic hexagonal lattice.
%%%%%%%%%%%%%%%%%%%%%%%%%%%%%%%%%%%%%%%%%%%%%%%%%%%%%%%%%%%%%%%%%%%%%%%%%%%%%%%%%%%%%%%%%%%%%%%%%%%%
%%%%%%%%%%%%%%%%%%%%%%%%%%%%%%%%%%%%%%%%%%%%%%%%%%%%%%%%%%%%%%%%%%%%%%%%%%%%%%%%%%%%%%%%%%%%%%%%%%%%
Its fundamental domain contains two vertices and it has three orbits of edges with labels $a,b$, and $c$ and potential values $u$ and $v$
as shown.

These are parameters for the tight binding model~\cite{AshcroftMermin,Kittel}, which involves the discrete Schr\"odinger
operator acting on functions $f\colon\calV\to\CC$. 
For a labeling $p\in(\RR^\times)^E\times\RR^W$, the corresponding \demph{discrete Schr\"odinger operator} \defcolor{$H_p$} acts on 
$f\colon \calV \rightarrow \CC$ as follows: 
 \begin{equation}
  \label{Eq:DSO}
   \bigl(H_p(f)\bigr)(u) \ = \  V(u) f(u)\ +\ \sum_{u\sim v} e(u,v) \cdot f(v)\,.
 \end{equation}
Here the sum is over all vertices $v \in \calV$ for which there is an edge between $u$ and $v$, written
\defcolor{$u{\sim}v$}.
The operator $H_p$ is a potential plus a weighted graph Laplacian.
In the standard formulation, the Laplacian contributes a sum of terms $-e(u,v)f(u)$ to~\eqref{Eq:DSO}.
In the formulation~\eqref{Eq:DSO}, that sum $-\sum_v e(u,v)$ is absorbed into the potential $V(u)$.
This gives a self-adjoint operator on the Hilbert space \defcolor{$\ell_2(\Gamma)$} of square summable functions
on $\calV$.
Its spectrum \defcolor{$\sigma(H_p)$} is a union of intervals in $\RR$.

Important features of the spectrum of $H_p$ are revealed by the Floquet transform $f\mapsto\hat{f}$, which is a version of the
Fourier transform.
References for Floquet theory include~\cite{BerkKuch, KuchBook, AAPGO}.  
Let $\defcolor{\TT}\subset\CC^\times$ be the group of unit complex numbers, which we identify with the set of unitary characters on $\ZZ$
via $z(n)=z^n$ for $n\in\ZZ$ and $z\in\TT$.
Similarly, $\TT^d$ is the group of unitary characters on $\ZZ^d$.
For $\alpha\in\ZZ^d$ and $z\in\TT^d$, we write \defcolor{$z^\alpha$} for $z_1^{\alpha_1}\dotsb z_d^{\alpha_d}$.
The Floquet transform of a function $f\in \ell_2(\Gamma)$ is a function \defcolor{$\hat{f}$} on $\TT^d\times\calV$
that is quasiperiodic in that 
$\hat{f}(z,u+\alpha) = z^\alpha \cdot \hat{f}(z,u)$ for $u\in\calV$ and $\alpha \in \ZZ^d$.
Such a quasiperiodic function is determined by its values on $\TT^d\times W$, which identifies the space of Floquet transforms with 
$L^2(\TT^d,\CC^W)=L^2(\TT^d)^W$, the space of square
integrable functions on $\TT^d$ with values in $\CC^W$. 
The Floquet transform is a linear isometry between $\ell_2(\Gamma)$ and $L^2(\TT^d)^W$.

As the labeling is $\ZZ^d$-periodic, $H_p$ commutes with the $\ZZ^d$-action and hence with the Floquet transform.
For $u\in W$, we have 
 \begin{equation}
  \label{Eq:L_Floquet}
   \bigl(H_p(\hat{f})\bigr)(z,u) \ =\
      V(u) \hat{f}(z,u) + \sum_{u{\sim}v+\alpha} e(u,v{+}\alpha) \cdot   z^\alpha \cdot \hat{f}(z,v)\,,
 \end{equation}
where the sum is over $v\in W$ and $\alpha\in\ZZ^d$ such that $u$ is adjacent to $v+\alpha$.
By the expression~\eqref{Eq:L_Floquet}, the operator $H_p$ is represented by a $|W| \times |W|$
\demph{Floquet matrix}, \defcolor{$H_p(z)$}, whose entries are Laurent polynomials in $z$. 
As $u\sim v{+}\alpha$ if and only if $v\sim u{-}\alpha$ and both edges have the same label, the Floquet matrix satisfies the identity
 \begin{equation} \label{eq:hermitianOnS1}
   H_p(z)^T\ =\ H_p(z^{-1})\,.
 \end{equation}
 Since $z\in\TT^d$, we have $\overline{z}=z^{-1}$, so that the Floquet matrix $H_p(z)$ is Hermitian.

%%%%%%%%%%%%%%%%%%%%%%%%%%%%%%%%%%%%%%%%%%%%%%%%%%%%%%%%%%%%%%%%%%%%%%%%%%%%%%%%%%%%%%%%%%%%%%%%%%%%
\begin{Example} \label{ex:floquethexagonal} 
Writing $z=(x,y)$ for points in $\TT^2$, the Floquet matrix for the hexagonal lattice of Figure~\ref{fig:HexGraph} is
 \begin{equation}
     \label{Eq:Matrix_For_Hexagonal_Lattice}
    H_p(x,y)\ =\ \left( \begin{matrix}
              u       & a + b x^{-1} + c y^{-1}\\
       a + b x + c y  &       v 
              \end{matrix}\right)\ . \qedhere
 \end{equation}
\end{Example}
%%%%%%%%%%%%%%%%%%%%%%%%%%%%%%%%%%%%%%%%%%%%%%%%%%%%%%%%%%%%%%%%%%%%%%%%%%%%%%%%%%%%%%%%%%%%%%%%%%%%

We introduce the geometric objects we want to study.
Write \defcolor{$D_p(z,\lambda)$} for the characteristic polynomial $\det(\lambda\cdot{\rm id}_{|W|} -H_p(z))$, which is called the
\demph{dispersion  polynomial}.
 
%%%%%%%%%%%%%%%%%%%%%%%%%%%%%%%%%%%%%%%%%%%%%%%%%%%%%%%%%%%%%%%%%%%%%%%%%%%%%%%%%%%%%%%%%%%%%%%%%%%%
\begin{Definition}[Real Bloch variety] \label{def:realbloch}
The \demph{real Bloch variety} of $H_p$ is the hypersurface 
\[
\calB \ =\ \calB_p \   \vcentcolon=\ \{ (z, \lambda) \in \TT^d \times \RR \mid D_p(z,\lambda) = 0 \}\,.
\]
A level set $\calF_p(\lambda) = \{ z\in\TT^d\mid  D_p(z,\lambda) = 0 \}$ of the Bloch variety is a
\demph{real Fermi~variety}.
\end{Definition}
%%%%%%%%%%%%%%%%%%%%%%%%%%%%%%%%%%%%%%%%%%%%%%%%%%%%%%%%%%%%%%%%%%%%%%%%%%%%%%%%%%%%%%%%%%%%%%%%%%%%

A reason this is important is that the \demph{spectrum $\sigma(H_p)$} of the Schr\"odinger operator $H_p$ is the projection of the real Bloch
variety to the $\RR$-axis.
Thus the Bloch variety encodes the interaction of the spectrum of $H_p$ with the action of $\ZZ^d$ on $\ell_2(\Gamma)$.

%%%%%%%%%%%%%%%%%%%%%%%%%%%%%%%%%%%%%%%%%%%%%%%%%%%%%%%%%%%%%%%%%%%%%%%%%%%%%%%%%%%%%%%%%%%%%%%%%%%%
\begin{Example}
Figure~\ref{FigTwo} shows two Bloch varieties for the hexagonal lattice with potential $V=0$, along with the spectrum of the
operator.
For display purposes, we represent $\TT^2$ by a fundamental domain $[-\frac{\pi}{2},\frac{3\pi}{2}]^2$ in $\RR^2$.
The edge parameters are $(a,b,c) = (6, 3, 2)$ on the left and on the right they are $(a,b,c) = (1, 1, 1)$, giving the graph Laplacian.
Its two singular points are called Dirac points.
%%%%%%%%%%%%%%%%%%%%%%%%%%%%%%%%%%%%%%%%%%%%%%%%%%%%%%%%%%%%%%%%%%%%%%%%%%%%%%%%%%%%%%%%%%%%%%%%%%%%
\begin{figure}[htb]
  \centering
   \begin{picture}(183,160)(-35,-10)
     \put(-2,-10){\includegraphics[height=159pt]{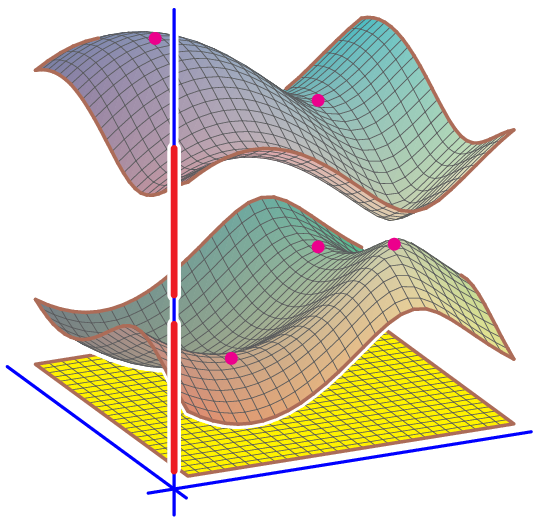}}
     \put(141,7){\small$x$} \put(5,21){\small$y$}
     \put(55,142){\small$\RR$}    \put(155,25){\small$\TT^2$}
     \put(-36,47){\small$\sigma(H)$}
     \thicklines
       \put(-10,46){{\color{white}\line(3,-1){56}}}
       \put(-10,45.5){{\color{white}\line(3,-1){56}}}
       \put(-10,45){{\color{white}\line(3,-1){56}}}
       \put(-10,44.5){{\color{white}\line(3,-1){56}}}
       \put(-10,44){{\color{white}\line(3,-1){56}}}
       
       \put(-10,55.5){{\color{white}\vector(3,1){44}}}
       \put(-10,55){{\color{white}\vector(3,1){44}}}
       \put(-10,54.5){{\color{white}\vector(3,1){44}}}
     \thinlines
       \put(48,25.666){\color{white}{\circle*{1.5}}}
       \put(47,26){\color{white}{\circle*{2}}}
       \put(45.5,26.1666){\color{white}{\circle*{2}}}
       \put(45,26.833){\color{white}{\circle*{3.5}}}
       \put(-10,55){\vector(3,1){58}}    \put(-10,45){\vector(3,-1){58}}
\end{picture}
 \qquad
\begin{picture}(190,160)(-19,-10)
    \put(0,1){\includegraphics[height=145pt]{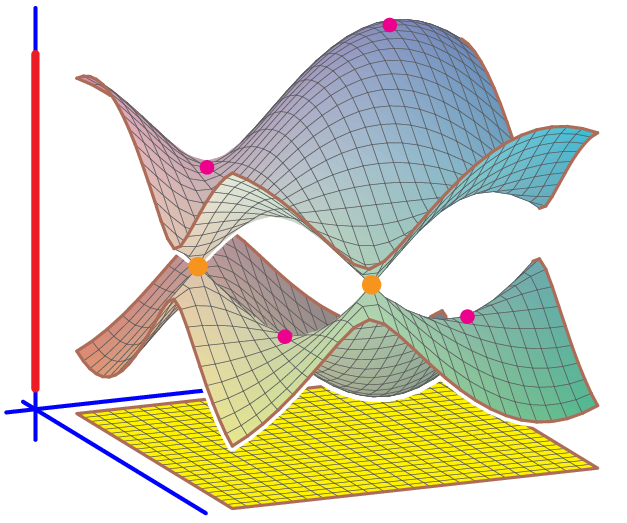}}
    \put(-19,81){\small$\sigma(H)$}
    \put(14,137){\small$\RR$}  
\end{picture}
\caption{Bloch varieties with edge parameters (6,3,2) and (1,1,1).}
\label{FigTwo}
\end{figure}
\end{Example}
%%%%%%%%%%%%%%%%%%%%%%%%%%%%%%%%%%%%%%%%%%%%%%%%%%%%%%%%%%%%%%%%%%%%%%%%%%%%%%%%%%%%%%%%%%%%%%%%%%%%

The terminology \demph{real} in Definition~\ref{def:realbloch} deserves some clarification.
Although $D_p(z, \lambda)$ has coefficients in $\RR$, it is customary in algebraic geometry to study the {\it complex}
points of the hypersurface defined by the equation $D_p(z,\lambda) = 0$.
Here is the corresponding definition. 

%%%%%%%%%%%%%%%%%%%%%%%%%%%%%%%%%%%%%%%%%%%%%%%%%%%%%%%%%%%%%%%%%%%%%%%%%%%%%%%%%%%%%%%%%%%%%%%%%%%%
\begin{Definition}
The \demph{(complex) Bloch variety} of $H_p$ is the hypersurface 
\[
\calB_{p}(\CC) \, = \, \{ (z, \lambda) \in (\CC^\times)^d \times \CC \mid D_p(z,\lambda) = 0 \}\,.
\]
The \demph{(complex) Fermi variety} $\calF_p(\lambda, \CC)$ is the level set of the Bloch variety at $\lambda\in\CC$.%\hfill$\diamond$
\end{Definition}
%%%%%%%%%%%%%%%%%%%%%%%%%%%%%%%%%%%%%%%%%%%%%%%%%%%%%%%%%%%%%%%%%%%%%%%%%%%%%%%%%%%%%%%%%%%%%%%%%%%%

Complexification has conceptual and computational advantages.
However, not all complex values of $z, \lambda$ are physically meaningful.
In applications, one is often interested in real solutions.
Abstractly, a point $(z, \lambda)$ is \demph{real} if it is fixed by a given antiholomorphic involution.
The most common involution is complex conjugation: $(z,\lambda) \mapsto (\overline{z}, \overline{\lambda})$.
As $D_p$ has real coefficients, $\calB_p(\CC)$ is stable under this involution, meaning that $(z,\lambda)$ belongs to $\calB_p(\CC)$ if
and only if $(\overline{z},\overline{\lambda})$ does.
The fixed points are the ordinary real points $\calB_p(\CC) \cap ((\RR^\times)^d \times \RR)$.
In spectral theory the meaningful choice of involution is $(z, \lambda) \mapsto (\overline{z^{-1}}, \overline{\lambda})$.
Since $D_p(z,\lambda) = D_p(z^{-1},\lambda)$, $\calB_p(\CC)$ is stable under this involution as well.
The real points corresponding to this involution are precisely the points in the real Bloch variety $\calB_p$. 

The choice of an antiholomorphic involution (real structure) on a complex algebraic variety $X$ determines which points of $X$ are real.
Most computational algorithms in real algebraic geometry assume the standard real structure corresponding to complex conjugation.
Bloch varieties motivate the need to develop tools to study natural non-standard real structures,
such as  $(z, \lambda) \mapsto (\overline{z^{-1}}, \overline{\lambda})$.

The coordinate projection $\lambda\colon(\CC^\times)^d\times \CC\to\CC$ restricts to a function on the complex Bloch variety $\calB_p(\CC)$.
Taking cues from polynomial optimization, in~\cite{FS24} it is shown that the \demph{critical points} of the function
$\lambda$ are those $(z,\lambda)$ which satisfy the critical point equations,
 \begin{equation}
  \label{eq:critpoints}
     D_p \ = \ z_1\frac{\partial D_p}{\partial z_1} \ = \ \cdots \ = \ z_d\frac{\partial D_p}{\partial z_d} \ = \ 0\,. 
 \end{equation}

A motivation to study critical points of $\lambda$ is the \demph{spectral edges nondegeneracy conjecture}~\cite[Conj.\ 5.25]{KuchBAMS},
which asserts that for generic values of the parameters  $p\in(\RR^\times)^E\times\RR^W$,
the extrema of the function $\lambda$ on each branch of the real Bloch variety are  nondegenerate.
Important notions in physics depend upon this assumption, and it has recently been proven for periodic operators on
$\RR^d$~\cite{TY2026}, but it remains open for discrete periodic operators.
Other geometric features, such as irreducibility of the Fermi variety~\cite{FL-G,FLM22,KuchmentVainberg2000} and the dimension
of the critical point locus~\cite{FaustKachkovskiy,FaustLiu,FS25}, are also of interest to mathematical physics.

We introduce some objects from geometric combinatorics related to the
dispersion polynomial $D_p(z,\lambda)$.
The exponents of monomials in $D_p(z,\lambda)$ form a subset of $\ZZ^d\times\NN$, called its \demph{support}.
Their convex hull is its \demph{Newton polytope}, $\calN = \defcolor{\calN(D_p(z,\lambda))}$.
This structure underlies the formulation of the critical point equations~\eqref{eq:critpoints}.
Multiplication by $z_i\in\CC^\times$
does not change the solution set in $(\CC^\times)^d\times\CC$ and the Newton polytope of
$z_i\frac{\partial D_p}{\partial z_i}$ is contained in the Newton polytope of $D_p$.
A result of~\cite{FS24} states that the number of critical points, counted with multiplicity,
is bounded above by the normalized volume of $\calN$, that is, its Euclidean volume
multiplied by $(d{+}1)!$. 

%%%%%%%%%%%%%%%%%%%%%%%%%%%%%%%%%%%%%%%%%%%%%%%%%%%%%%%%%%%%%%%%%%%%%%%%%%%%%%%%%%%%%%%%%%%%%%%%%%%%
\begin{Example} \label{ex:newtonhexagon}
Consider the dispersion polynomial $D_p = \det(\lambda \cdot {\rm id} - H_p)$, with $H_p$ as in \eqref{Eq:Matrix_For_Hexagonal_Lattice}.
Using a computer algebra system, such as \texttt{Macaulay2}~\cite{M2}, \texttt{Oscar}~\cite{Oscar}, or \texttt{Singular}~\cite{Singular}, 
we may check that for general values of the parameters $p$, the critical point equations~\eqref{eq:critpoints} have twelve complex solutions.
This equals the upper bound given in~\cite[Cor.~2.5]{FS24}, which is the normalized volume of $\calN$.
Figure~\ref{Fig:hexPolytope} shows the support of $D_p$ as the columns of a matrix and its Newton polytope, $\calN$.
%%%%%%%%%%%%%%%%%%%%%%%%%%%%%%%%%%%%%%%%%%%%%%%%%%%%%%%%%%%%%%%%%%%%%%%%%%%%%%%%%%%%%%%%%%%%%%%%%%%%
\begin{figure}[htb]
 \centering
 ${\displaystyle \left(\begin{array}{rrrrrrrrr}
     0&0&0& 1&0&-1&-1& 0& 1\\
     0&0&0& 0&1& 1& 0&-1&-1\\
     2&1&0& 0&0& 0& 0& 0& 0\end{array}\right)}$
   \qquad
   \raisebox{-25pt}{\begin{picture}(130,80)(-20,0)
       \put(-18,30){$\calN=$}
       \put(0,0){\includegraphics[height=60pt]{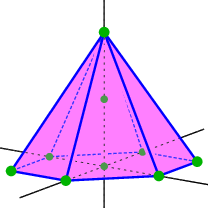}}
       \end{picture}}
   \caption{Support of dispersion polynomial of the hexagonal graph and its Newton polytope.}
   \label{Fig:hexPolytope}
\end{figure}
%%%%%%%%%%%%%%%%%%%%%%%%%%%%%%%%%%%%%%%%%%%%%%%%%%%%%%%%%%%%%%%%%%%%%%%%%%%%%%%%%%%%%%%%%%%%%%%%%%%%
The base of $\calN$ is a hexagon with area 3 and it is a pyramid of height 2.
Thus this pyramid has Euclidean volume 2 and hence normalized volume $3!\, \vol(\calN)=12$.
\end{Example}
%%%%%%%%%%%%%%%%%%%%%%%%%%%%%%%%%%%%%%%%%%%%%%%%%%%%%%%%%%%%%%%%%%%%%%%%%%%%%%%%%%%%%%%%%%%%%%%%%%%%

We are now ready to define our main object of study.

%%%%%%%%%%%%%%%%%%%%%%%%%%%%%%%%%%%%%%%%%%%%%%%%%%%%%%%%%%%%%%%%%%%%%%%%%%%%%%%%%%%%%%%%%%%%%%%%%%%%
\begin{Definition} \label{def:Bloch_Discriminant}
  Let $\Gamma$ be a $\ZZ^d$-periodic graph.
  A \demph{Bloch discriminant} for $\Gamma$ is a hypersurface
  $\BD_\Gamma \subset (\CC^\times)^E \times \CC^W$ such that the number of real solutions
  $(z,\lambda) \in \TT^d \times \RR$ of \eqref{eq:critpoints} is constant for all $p$ in each connected component  of
  $\left ( (\RR^\times)^E \times \RR^W\right) \smallsetminus \BD_\Gamma$.
  These connected components  are called \demph{discriminant chambers}.
\end{Definition}
%%%%%%%%%%%%%%%%%%%%%%%%%%%%%%%%%%%%%%%%%%%%%%%%%%%%%%%%%%%%%%%%%%%%%%%%%%%%%%%%%%%%%%%%%%%%%%%%%%%%

We shall write $\BD_\Gamma$ for a Bloch discriminant hypersurface in
$(\CC^\times)^E \times \CC^W$, as well as for a defining polynomial equation.
That is, the polynomial $\BD_\Gamma$ defines the Bloch discriminant variety
$\BD_\Gamma = \{ p \in (\CC^\times)^E \times \CC^W \mid \BD_\Gamma(p) = 0 \}$. 

We compute a first Bloch discriminant for the hexagonal lattice---summarized below in Theorem \ref{T:hexagonal} and expanded on in
Section \ref{S:Hexagonal}.

%%%%%%%%%%%%%%%%%%%%%%%%%%%%%%%%%%%%%%%%%%%%%%%%%%%%%%%%%%%%%%%%%%%%%%%%%%%%%%%%%%%%%%%%%%%%%%%%%%%%
\begin{theorem}\label{T:hexagonal}
  Let $\Gamma$ be the hexagonal lattice of Figure \ref{fig:HexGraph}.
  The quintic polynomial 
\[
  \BD_\Gamma \ = \ (u-v)(a+b+c)(a+b-c)(a-b+c)(a-b-c)
\]
is a Bloch discriminant for $\Gamma$.
\end{theorem}
%%%%%%%%%%%%%%%%%%%%%%%%%%%%%%%%%%%%%%%%%%%%%%%%%%%%%%%%%%%%%%%%%%%%%%%%%%%%%%%%%%%%%%%%%%%%%%%%%%%%

Let us write $\defcolor{\Delta}$ for the factors of the discriminant polynomial $\BD_\Gamma$  involving only $a,b,c$.
Its sign dictates the critical point structure.
It is worth noting that $\Delta<0$ if and only if $|a|,|b|,|c|$ form the sides of a triangle.
We display this in the $(b,c)$-plane where $a=1$.
\[
  \begin{picture}(125,72)
    \put(0,0){\includegraphics{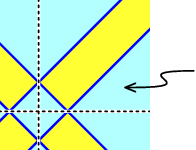}}
    \put(21,64){\small$c$}  \put(75,15){\small$b$}
    \put(95,35){$\RR^3_+$}
    \put(95,61){$\Delta<0$} \put(93,64){\vector(-1,0){31}}
   \end{picture}
\]
The dotted lines are where $b$ or $c$ vanish, the solid lines are where $\Delta=0$.
The discriminant $\Delta$ is negative in the yellow regions and positive in the light blue regions.

Table \ref{table:hexagonal} summarizes the landscape of the critical points based on corresponding parameter values.
This establishes the claim in Proposition~\ref{P:hexagonal_CPA} about the minimum of $\lambda^+$ when $\Delta<0$.
The triples in the Hessian column record the numbers of real nondegenerate critical points that are local minima, saddle points,
and local maxima.
Only the first two rows concern discriminant chambers; the remaining rows describe the landscape of critical points where
the Bloch discriminant vanishes.

%%%%%%%%%%%%%%%%%%%%%%%%%%%%%%%%%%%%%%%%%%%%%%%%%%%%%%%%%%%%%%%%%%%%%%%%%%%%%%%%%%%%%%%%%%%%%%%%%%%%
\begin{table}[htb]
\centering
\begin{tabular}{||c | c | c | c | c||} 
 \hline
 & $\Delta$ & nondegenerate   & degenerate      & Hessian \\ [0.3ex] 
 &          & critical points & critical points &  \\ [0.3ex] 
  \hline \hline 
  $u \neq v$ & $< 0$ & 12 & 0            & (3,6,3)\rule{0pt}{12pt} \\  [0.6ex]
             & $> 0$ &  8 & 0            & (2,4,2) \\ [0.6ex]
             & $= 0$ &  6 & 2 (mult.\ 3) & (1,4,1) \\ [0.6ex]
 \hline\hline 
 $u = v$ & $< 0$ & 8 & 2 (Dirac) & (1,6,1) \rule{0pt}{12pt}\\ [0.6ex]
         & $> 0$ & 8 & 0            & (2,4,2) \\ [0.6ex]
         & $= 0$ & 6 & 1 (mult.\ 6) & (1,4,1) \\ [0.2ex]
  \hline
\end{tabular}
\caption{Critical point landscape summary for the hexagonal lattice.\rule{0pt}{14pt}}
\label{table:hexagonal}
\end{table}
%%%%%%%%%%%%%%%%%%%%%%%%%%%%%%%%%%%%%%%%%%%%%%%%%%%%%%%%%%%%%%%%%%%%%%%%%%%%%%%%%%%%%%%%%%%%%%%%%%%%

Parameters such that the critical point equations~\eqref{eq:critpoints} have a real solution of multiplicity at least
  two lie in the Bloch Discriminant.
That is, the Jacobian matrix of the critical point equations will not have full rank at that solution.
We exploit this to give a procedure for computing a Bloch discriminant under the following assumption:
for generic parameter values, there are finitely many solutions to \eqref{eq:critpoints} and all solutions have
multiplicity one.
This assumption holds in the case studies from Sections \ref{S:Hexagonal}-\ref{S:five}.

%%%%%%%%%%%%%%%%%%%%%%%%%%%%%%%%%%%%%%%%%%%%%%%%%%%%%%%%%%%%%%%%%%%%%%%%%%%%%%%%%%%%%%%%%%%%%%%%%%%%
\begin{algorithm}[htb]
\caption{Computing a Bloch discriminant for $\Gamma$}
\label{alg:main}
\begin{algorithmic}[1]

\State Work in the polynomial ring $R = \QQ[z_1, \ldots, z_d, p_0, \ldots, p_m, \lambda, \mu]$ with $d + |E| + |W| + 2$ variables.
  Here, $p_i$ are the parameters labeling the edges and vertices of $\Gamma$.

\State Let $g_0, g_1, \ldots, g_d \in R$ be the numerators of the rational functions appearing in \eqref{eq:critpoints}.
  That is, these are the numerators of the dispersion polynomial and its toric derivatives.

\State Let $g_{d+1} = z_1z_2 \cdots z_d \, (\prod_{e \in E} p_e) \, \mu - 1 \in R$, where the product is over edge labels.

\State Compute the $(d+2) \times (d+2)$ Jacobian determinant of $g_0, g_1, \ldots, g_{d+1}$: 
\[ h \, = \, \det \begin{pmatrix}
\frac{\partial g_0}{\partial z_1} & \cdots & \frac{\partial g_0}{\partial z_d} & \frac{\partial g_0}{\partial \lambda} & \frac{\partial g_0}{\partial \mu} \\ 
\vdots & & \vdots & \vdots & \vdots \\
\frac{\partial g_d}{\partial z_1} & \cdots & \frac{\partial g_d}{\partial z_d} & \frac{\partial g_d}{\partial \lambda} & \frac{\partial g_d}{\partial \mu} \\ 
\frac{\partial g_{d+1}}{\partial z_1} & \cdots & \frac{\partial g_{d+1}}{\partial z_d} & \frac{\partial g_{d+1}}{\partial \lambda} & \frac{\partial g_{d+1}}{\partial \mu} \\ 
\end{pmatrix}\]

\State Let $I \subset R$ be the ideal $\langle g_0, g_1, \ldots, g_{d+1}, h \rangle$. 

\State Compute the elimination ideal $J = I \cap \QQ[p_0, \ldots, p_m]$. 

\State Return a nonzero generator of the radical ideal $\sqrt{J}$. 
\end{algorithmic}
\end{algorithm}
%%%%%%%%%%%%%%%%%%%%%%%%%%%%%%%%%%%%%%%%%%%%%%%%%%%%%%%%%%%%%%%%%%%%%%%%%%%%%%%%%%%%%%%%%%%%%%%%%%%%
\begin{proof}[Proof of correctness]
  Computer algebra systems work with ordinary polynomials.
  This necessitates taking numerators in Step 2 (the denominators are  monomials in $z$).
  The auxiliary indeterminate $\mu$ and the polynomial $g_{d+1}$ in Step 3 ensure that neither any coordinate $z_i$ nor any
  parameter $p_e, e \in E$  vanishes at a solution.
  The ideal $I$ vanishes at $(z,p,\lambda,\mu)$ which are singular solutions to the critical point equations with
  $z\in(\CC^\times)^d$ and $p\in(\CC^\times)^E \times \CC^W$.

  The elimination ideal $J$ of Step 6 defines the parameters $p$ for which the critical point equations have a singular
  solution.
  By assumption, this does not happen for generic parameters, so $J$ is not the zero ideal.
  When $J$ is principal, $J=\langle f\rangle$ and its radical is generated by the square-free part of $f$. If $J$ is not
  principal, then the locus of singular parameters has potentially embedded components of codimension $\geq 2$.
  In this case, the algorithm returns the defining polynomial of a hypersurface which contains that locus potentially
  strictly.
  This is a Bloch discriminant by definition. 
\end{proof}
%%%%%%%%%%%%%%%%%%%%%%%%%%%%%%%%%%%%%%%%%%%%%%%%%%%%%%%%%%%%%%%%%%%%%%%%%%%%%%%%%%%%%%%%%%%%%%%%%%%%

Applying Algorithm~\ref{alg:main} to the hexagonal lattice from Figure \ref{fig:HexGraph}, gives a polynomial of degree $13$ with nine
irreducible factors.
Five linear factors appear in Theorem \ref{T:hexagonal}.
The remaining four quadratic factors are 
\[
  (u-v)^2\ +\  4(a \pm b \pm c)^2\,.
\]
These polynomials will appear again in our analysis in Section \ref{S:Hexagonal}.
Notice that their zero loci in $(\mathbb{R}^\times)^E \times \mathbb{R}^W$ are contained in the hypersurface from Theorem
\ref{T:hexagonal}.
In that sense, these four factors are redundant. Theorem \ref{T:hexagonal} provides a more efficient Bloch~discriminant.

%%%%%%%%%%%%%%%%%%%%%%%%%%%%%%%%%%%%%%%%%%%%%%%%%%%%%%%%%%%%%%%%%%%%%%%%%%%%%%%%%%%%%%%%%%%%%%%%%%%%%
  The elimination in Step 6 of Algorithm~\ref{alg:main} may be prohibitively resource-intensive and thus Algorithm~\ref{alg:main} is only
  practical to compute small examples.
  We will exploit structures of the dispersion polynomial to develop more specialized tools, and we present several case studies.
  We additionally make an important observation about Fermi varieties.
  A critical value (the $\lambda$-coordinate of a critical point) is a \demph{critical energy}. 
%%%%%%%%%%%%%%%%%%%%%%%%%%%%%%%%%%%%%%%%%%%%%%%%%%%%%%%%%%%%%%%%%%%%%%%%%%%%%%%%%%%%%%%%%%%%%%%%%%%%

%%%%%%%%%%%%%%%%%%%%%%%%%%%%%%%%%%%%%%%%%%%%%%%%%%%%%%%%%%%%%%%%%%%%%%%%%%%%%%%%%%%%%%%%%%%%%%%%%%%%
\begin{lemma}\label{L:singularFermi}
  The Fermi variety  $\calF_p(\lambda, \CC)$ is singular if and only if $\lambda$ is a critical energy,
  and its singular points are the critical points of the  Bloch variety at the energy $\lambda$. 
\end{lemma}
%%%%%%%%%%%%%%%%%%%%%%%%%%%%%%%%%%%%%%%%%%%%%%%%%%%%%%%%%%%%%%%%%%%%%%%%%%%%%%%%%%%%%%%%%%%%%%%%%%%%
\begin{proof}
  The Fermi variety
  $\calF_p(\lambda_0, \CC)$ is singular at~$z_0$ if and only if $(z_0,\lambda_0)$ satisfies \eqref{eq:critpoints}.
\end{proof}
%%%%%%%%%%%%%%%%%%%%%%%%%%%%%%%%%%%%%%%%%%%%%%%%%%%%%%%%%%%%%%%%%%%%%%%%%%%%%%%%%%%%%%%%%%%%%%%%%%%%

As studied in~\cite{GKZ}, a hypersurface is singular if and only if an appropriate discriminant vanishes at its coefficients.
This naturally leads to the topic of our next section.

%%%%%%%%%%%%%%%%%%%%%%%%%%%%%%%%%%%%%%%%%%%%%%%%%%%%%%%%%%%%%%%%%%%%%%%%%%%%%%%%%%%%%%%%%%%%%%%%%%%%
\section{Discriminants of reciprocal polynomials} \label{S:ReciprocalDiscriminants}

A Laurent polynomial $f(z)$ is  \demph{reciprocal} if it is invariant under inversion of indeterminates, $f(z) = f(z^{-1})$.
As a polynomial in $z$, the dispersion polynomial $D_p(z,\lambda)$ is reciprocal.
The hypersurface $\calV(f)\subset(\CC^\times)^d$ of a reciprocal polynomial is invariant under the map $z\mapsto z^{-1}$.
Let us consider the family of  hypersurfaces defined by the reciprocal polynomials
\begin{equation} \label{eq:fzp}
  f(z,p) \ \vcentcolon= \ p_0 + p_1 \, \left ( \frac{z^{a_1} + z^{-a_1}}{2} \right) + \cdots + p_m \,
  \left ( \frac{z^{a_m} + z^{-a_m}}{2} \right )\ .
\end{equation}
Here, the exponents $\{a_1, \ldots, a_m\}\subset\ZZ^d$ are nonzero and distinct and we assume that they span $\RR^d$.
We consider $p$ to be a point in the projective space $\PP^m$.
Let $A\vcentcolon=\{0,\pm a_1,\dotsc,\pm a_m\}$ be the exponents appearing in~\eqref{eq:fzp}.
The \demph{discriminant variety} $\defcolor{\nabla_r(A)}\subset\PP^m$ is the closure of the set of points $p\in\PP^m$ such that
the hypersurface $\calV(f(z,p))\subset(\CC^\times)^d$ is singular.
This is a reducible hypersurface, as shown in \cite{RDisc}.
Its primitive squarefree defining polynomial is denoted by $\defcolor{\Delta_r(A)}\in\ZZ[p]$, which is well-defined up to
multiplication by a unit $\pm1$ of $\ZZ$. 
We offer a method to compute a polynomial multiple of $\Delta_r(A)$ in Algorithm~\ref{alg:discr}.

%%%%%%%%%%%%%%%%%%%%%%%%%%%%%%%%%%%%%%%%%%%%%%%%%%%%%%%%%%%%%%%%%%%%%%%%%%%%%%%%%%%%%%%%%%%%%%%%%%%%
\begin{algorithm}[htb]
\caption{Computing discriminants}
\label{alg:discr}
\begin{algorithmic}[1]

\State Let $R = \QQ[z_1, \ldots, z_d,  w_1,\dotsc,w_d,p_0, \ldots, p_m]$.

\State Let $I$ be the ideal whose generators are the (numerators of) the polynomial $f(z,p)$ and its partial derivatives with respect to each
$z_i$, as well as the polynomials $z_i w_i - 1$.

\State Compute a Gr\"obner basis $G$ for $I$ using the degree reverse lexicographic term order.

\State If there is a nonzero polynomial in $G\cap \QQ[p]$ return its squarefree part and stop.

\State Otherwise, compute the elimination ideal $J=I\cap \QQ[p]$ and return the squarefree part of its generator.
\end{algorithmic}
\end{algorithm}
%%%%%%%%%%%%%%%%%%%%%%%%%%%%%%%%%%%%%%%%%%%%%%%%%%%%%%%%%%%%%%%%%%%%%%%%%%%%%%%%%%%%%%%%%%%%%%%%%%%%

The discriminant $\Delta_r(A)$ will be the squarefree part of the polynomial of smallest degree in
$I\cap\QQ[p_0,\dotsc,p_m]$ and hence Step 5 will compute it.
It is not clear that  a multiple of it
occurs in the Gr\"obner basis computed in Step 3, in which case Algorithm~\ref{alg:discr} stops at Step 4 with a polynomial
  multiple of $\Delta_r(A)$.
We used this to compute discriminants in this paper. Remarkably, in every case, this stops at Step 4, returning $\Delta_r(A)$.

This discriminant variety $\nabla_r(A)$ is a reducible hypersurface:
The reciprocal structure causes singular points at the fixed points of $z \mapsto z^{-1}$ which lie on $\calV(f(z,p))$.
This is because each partial derivative  $\frac{\partial f}{\partial z_j}$ vanishes at such a point.
Let us denote this set of $2^d$ fixed points by $\defcolor{\CP} \vcentcolon= \{ (\pm 1, \dotsc, \pm 1)\} \subset \TT^d$,
which are called \demph{corner points} in~\cite{BCCM_LocalGlobal}. 

%%%%%%%%%%%%%%%%%%%%%%%%%%%%%%%%%%%%%%%%%%%%%%%%%%%%%%%%%%%%%%%%%%%%%%%%%%%%%%%%%%%%%%%%%%%%%%%%%%%%
\begin{lemma} \label{lem:linfacdeltarec}
  If any of the linear forms $\{f(z^*,p) \mid  z^* \in \CP \}$  vanishes at the point $p \in \PP^m$, then the
  hypersurface $\calV(f(z,p)) \subset (\mathbb{C}^\times)^d$ is singular at $z^*$. 
\end{lemma}
%%%%%%%%%%%%%%%%%%%%%%%%%%%%%%%%%%%%%%%%%%%%%%%%%%%%%%%%%%%%%%%%%%%%%%%%%%%%%%%%%%%%%%%%%%%%%%%%%%%%

Thus each hyperplane $\calV(f(z^*,p)) \subset \mathbb{P}^{m}$ for $z^*\in \CP$ is an irreducible component of $\nabla_r(A)$.
By Lemma \ref{lem:linfacdeltarec}, $\Delta_r(A)$ is reducible and has at least the irreducible factors $\{ f(z^*,p)\mid z^*\in\CP\}$.
Another contribution is captured by the \demph{reciprocal discriminant}.
Let $\defcolor{Z_A}\subset(\CC^\times)^d$ be the hypersurface defined by
$(z^{a_1} - z^{-a_1})(z^{a_2} - z^{-a_2}) \cdots (z^{a_m} - z^{-a_m}) = 0$, which is a torus arrangement containing $\CP$.

%%%%%%%%%%%%%%%%%%%%%%%%%%%%%%%%%%%%%%%%%%%%%%%%%%%%%%%%%%%%%%%%%%%%%%%%%%%%%%%%%%%%%%%%%%%%%%%%%%%%
\begin{Definition} \label{def:recdisc}
  The \demph{reciprocal $A$-discriminant variety} \defcolor{$\nabla_A^\rec$} is the Zariski closure of the set
\[
   \Bigl\{ p \in \PP^m  \mid  f = z_1\frac{\partial f}{\partial z_1} = \dotsb
     = z_d\frac{\partial f}{\partial z_d} = 0, \text{ for some } z \in (\CC^\times)^d \setminus Z_A \Bigr\}\,,
\]
 where $f=f(z,p)$ is the reciprocal polynomial~\eqref{eq:fzp}.  
\end{Definition}
%%%%%%%%%%%%%%%%%%%%%%%%%%%%%%%%%%%%%%%%%%%%%%%%%%%%%%%%%%%%%%%%%%%%%%%%%%%%%%%%%%%%%%%%%%%%%%%%%%%%

For generic $p\in\nabla_A^\rec$, the singular points of $\calV(f(z,p))$ do not occur at a corner point.
They also do not occur at points $z$ of the torus arrangement $Z_A$.

%%%%%%%%%%%%%%%%%%%%%%%%%%%%%%%%%%%%%%%%%%%%%%%%%%%%%%%%%%%%%%%%%%%%%%%%%%%%%%%%%%%%%%%%%%%%%%%%%%%%
\begin{Remark} \label{rem:gkz}
  The classical theory of \demph{$A$-discriminants}  developed by Gel'fand, Kapranov and Zelevinsky~\cite{GKZ}
  also associates a discriminant to a set $A=\{a_0,\dotsc,a_n\}$ of integer exponent vectors, which vanishes when a
  corresponding hypersurface is singular.
  Specifically, let   $g(z,q) \, = \, q_0 \, z^{a_0} + q_1 \, z^{a_1} + \cdots + q_n \, z^{a_n}$ be a family of Laurent polynomials.
  The \demph{$A$-discriminant variety} $\defcolor{\nabla_A}\subset\PP^n$ is the Zariski closure of the set
  \[
    \bigl\{ q \in \PP^n \mid g(z,q) = \frac{\partial g}{\partial z_1}(z,q) = \cdots = \frac{\partial g}{\partial z_d}(z,q) = 0,
    \text{ for some } z \in (\CC^\times)^d \bigr\}\,.
  \]
  When this is a hypersurface, its defining polynomial is the $A$-discriminant, $\defcolor{\Delta_A}\in\ZZ[q]$.
  It is also identified with the projective dual variety to the toric variety parameterized by monomials in $z$ with exponents in $A$.
  We will illustrate this in Example \ref{ex:reciprocal_disc_hexagonal} below. 
\end{Remark}
%%%%%%%%%%%%%%%%%%%%%%%%%%%%%%%%%%%%%%%%%%%%%%%%%%%%%%%%%%%%%%%%%%%%%%%%%%%%%%%%%%%%%%%%%%%%%%%%%%%%

%%%%%%%%%%%%%%%%%%%%%%%%%%%%%%%%%%%%%%%%%%%%%%%%%%%%%%%%%%%%%%%%%%%%%%%%%%%%%%%%%%%%%%%%%%%%%%%%%%%%
\begin{proposition} \label{prop:irreducible}
  The reciprocal $A$-discriminant variety $\nabla_A^\rec \subset \PP^m$ is irreducible of codimension at least one. 
\end{proposition}
%%%%%%%%%%%%%%%%%%%%%%%%%%%%%%%%%%%%%%%%%%%%%%%%%%%%%%%%%%%%%%%%%%%%%%%%%%%%%%%%%%%%%%%%%%%%%%%%%%%%
\begin{proof}
  Let $f=f(z,p)$ be the reciprocal polynomial~\eqref{eq:fzp}.
  Consider the incidence variety 
\[
  \defcolor{Y} \ \vcentcolon=\ 
   \bigl\{ (z,p) \in ((\CC^\times)^d \setminus Z_A) \times \PP^m \mid
   f = z_1\frac{\partial f}{\partial z_1} = \dotsb = z_d\frac{\partial f}{\partial z_d} = 0 \bigr\}\,.
\]
The variety $\nabla_A^\rec$ is the closure of the image of the projection $Y \rightarrow \PP^m$.
Hence it suffices to show that $Y$ is irreducible of dimension $\leq m-1$.
For this, consider the other projection $\pi \colon Y \rightarrow (\CC^\times)^d \smallsetminus Z_A$.
A fiber $\pi^{-1}(z)$ of that projection is a linear space in $\PP^m$ as the equations defining $Y$ are linear in $p$.

Inspired by~\cite{MFRCSD}, we write the equations for $Y$ propitiously.
Let $B$ be the $d\times(m{+}1)$-matrix with columns $0, a_1,\dotsc,a_m$.
Let $\defcolor{\tau(z)}\vcentcolon= (2,z^{a_1}+z^{-a_1},\dotsc,z^{a_m}+z^{-a_m})$ be the coefficients of $p$ in $2f(z,p)$, and let
$\defcolor{\tau'(z)}\vcentcolon= (0,z^{a_1}-z^{-a_1},\dotsc,z^{a_m}-z^{-a_m})$.
Writing $\tau(z)$ as a row vector and the coefficients $p$ as a column vector,
we have
 \begin{equation}\label{Eq:vectorised_f}
   2f(z,p)\ =\ \tau(z) \cdot p\,.
 \end{equation}
The vector $(2z_i\partial f/\partial z_i\mid i=1,\dotsc,d)$ of toric derivatives is 
 \begin{equation}\label{Eq:vectorised_derivs}
   B \cdot \mbox{{\rm Diag}}(\tau'(z)) \cdot p\,.
 \end{equation}
For a fixed $z\in(\CC^\times)^d$, the $1{+}d$ linear equations $\tau(z)\cdot p=0$ and  $B \cdot \mbox{{\rm Diag}}(\tau'(z)) \cdot p=0$
define $\pi^{-1}(z)$.
If $z\not\in Z_A$, then this set of linear equations has rank $d+1$.
Indeed, the coefficient of $p_0$ is 2 in the first  and 0 in the remaining equations.
Our assumption on the exponents is that the matrix $B$ has rank $d$, and as $z\not\in Z_A$, no entry of $\mbox{{\rm Diag}}(\tau'(z))$ vanishes
except for the initial entry corresponding to $p_0$.
This completes the proof.
\end{proof}
%%%%%%%%%%%%%%%%%%%%%%%%%%%%%%%%%%%%%%%%%%%%%%%%%%%%%%%%%%%%%%%%%%%%%%%%%%%%%%%%%%%%%%%%%%%%%%%%%%%%

If $X\subset\PP^m$ is the variety parameterized by 
 \[ 
    (\CC^\times)^d\ \ni\ z \  \longmapsto \  \left ( 1\;:\;\frac{z^{a_1} + z^{-a_1}}{2}\;:\; \dotsb \;:\; \frac{z^{a_m} + z^{-a_m}}{2} \right )\,,
 \]
then $\nabla_A^{\rm rec}$ is (the closure of) the set of hyperplanes in $\PP^m$ that are tangent to a smooth point of $X$.
That is, it is the projective dual to $X$.
Bertini's Theorem implies that it is a proper subvariety.
Note that $X$ is the projective closure of the affine variety parameterized by
$u \mapsto (\cos(a_1 \cdot u), \ldots, \cos(a_m \cdot u))$, which was called a \demph{Chebyshev variety} in \cite{bel2024chebyshev}. 

Proposition \ref{prop:irreducible} justifies the following definition. 

%%%%%%%%%%%%%%%%%%%%%%%%%%%%%%%%%%%%%%%%%%%%%%%%%%%%%%%%%%%%%%%%%%%%%%%%%%%%%%%%%%%%%%%%%%%%%%%%%%%%
\begin{Definition}
  If $\nabla_A^\rec$ is a hypersurface, then the \demph{reciprocal $A$-discriminant} is the irreducible polynomial
  $\defcolor{\Delta_A^\rec} \in \ZZ[p_0, \ldots, p_m]$  whose complex variety is 
  $\nabla_A^\rec$.
  %
  %  Irreducibility implies that the coefficeints do not have a common nonunit factor in $\ZZ$.
  %
  This polynomial is well-defined up to multiplication by $-1$.
  If $\nabla_A^\rec$ is not a hypersurface, then we set $\Delta_A^\rec = 1$. 
\end{Definition}
%%%%%%%%%%%%%%%%%%%%%%%%%%%%%%%%%%%%%%%%%%%%%%%%%%%%%%%%%%%%%%%%%%%%%%%%%%%%%%%%%%%%%%%%%%%%%%%%%%%%

We observe that the reciprocal discriminant variety $\nabla_A^\rec$ need not be a hypersurface.

%%%%%%%%%%%%%%%%%%%%%%%%%%%%%%%%%%%%%%%%%%%%%%%%%%%%%%%%%%%%%%%%%%%%%%%%%%%%%%%%%%%%%%%%%%%%%%%%%%%%
\begin{Example}
  Consider the following reciprocal polynomial in three variables:
  \[
  f(z,p)\ \vcentcolon=\  p_0 + p_1 \, \frac{z_1 + z_1^{-1}}{2} + p_2 \, \frac{z_2 + z_2^{-1}}{2}
  + p_3 \, \frac{z_1z_2^{-1} + z_1^{-1}z_2}{2} + p_4 \, \frac{z_3 + z_3^{-1}}{2}\,.
  \]
  Its set of exponents form the columns of the matrix
  \begin{equation}\label{Eq:3x9}
    A\ \vcentcolon=\ \left(\begin{array}{rrrrrrrrr}
      0 & 1 & -1 & 0 &  0 &  1 & -1 & 0 &  0\\
      0 & 0 &  0 & 1 & -1 & -1 &  1 & 0 &  0\\
      0 & 0 &  0 & 0 &  0 &  0 &  0 & 1 & -1 \end{array}\right)\ .
  \end{equation}
  We use Algorithm~\ref{alg:discr} to compute the discriminant $\Delta_r(A)\subset\PP^4$, which
  is a hypersurface with $10$ components.
  Eight are the hyperplanes of Lemma~\ref{lem:linfacdeltarec}.
  The other two are quartic threefolds $C_+$ and $C_-$ defined by 
  \[
  2 p_0 p_1 p_2 p_3\ -\ p_1^2 p_2^2\  -\ p_1^2 p_3^2\  -\ p_2^2 p_3^2\  \pm\ 2 p_1 p_2 p_3 p_4 \ = \ 0\,.
  \]
  The parameters in $C_\pm$ give surfaces whose singular locus lies in the translated subtori $z_3=\pm1$.
  
  The reciprocal $A$-discriminant has codimension two.
  It is an irreducible component of the intersection of these quartics, defined by
  \begin{equation}\label{Eq:rec_A-discrim}
    2 p_0 p_1 p_2 p_3\ -\ p_1^2 p_2^2\ -\ p_1^2 p_3^2\  -\ p_2^2 p_3^2\  =\  p_4\ =\ 0\,.
  \end{equation}
  As $p_4=0$, the singular polynomial $f(z,p)$ depends only upon $z_1$ and $z_2$.
  A hypersurface parameterized by a point $p$ in the reciprocal $A$-discriminant has singular locus a union of 
  lines $(z_1^*,z_2^*)\times\CC^\times$ where $(z_1^*,z_2^*)$ is a singular point of the hypersurface in $(\CC^\times)^2$ defined by
  $f(z,p)$ (as $p_4=0$, the indeterminate $z_3$ does not appear in $f(z,p)$).
\end{Example}
%%%%%%%%%%%%%%%%%%%%%%%%%%%%%%%%%%%%%%%%%%%%%%%%%%%%%%%%%%%%%%%%%%%%%%%%%%%%%%%%%%%%%%%%%%%%%%%%%%%%

This example shows that $\nabla_A^\rec$ may not be a hypersurface and that the discriminant $\Delta_r(A)$ may have
other irreducible components than given by Lemma~\ref{lem:linfacdeltarec} or Definition~\ref{def:recdisc}.\smallskip

We now study $\Delta_r(A)$ where the exponents $A$ appear in dispersion polynomials for certain graphs,
beginning with the hexagonal lattice. 

%%%%%%%%%%%%%%%%%%%%%%%%%%%%%%%%%%%%%%%%%%%%%%%%%%%%%%%%%%%%%%%%%%%%%%%%%%%%%%%%%%%%%%%%%%%%%%%%%%%%
\begin{Example} \label{ex:reciprocal_disc_hexagonal}
  The dispersion polynomial $D_p = \det (\lambda\cdot{\rm id} - H_p)$, with $H_p$ the Floquet matrix of Example~\ref{ex:floquethexagonal},
  belongs to the following family of Laurent polynomials, 
 \[
     g(z,q) \ = \  q_0\ +\ q_1 \, x\ +\ q_2 \, x^{-1} \ +\ q_3 \, y\ +\ q_4\, y^{-1}\ +\ q_5\,xy^{-1}\  +\ q_6\, x^{-1}y\,.
 \]
  The exponents appearing in this polynomial are the columns of the matrix 
\[
  A \ = \ \left(\begin{array}{rrrrrrr}
0 & 1 & -1 & 0 &  0 &  1 & -1  \\ 
0 & 0 &  0 & 1 & -1 & -1 &  1 
  \end{array}\right)\,.
\]
These are the lattice points in the polygon which is the base of the pyramid shown in Figure~\ref{Fig:hexPolytope}, and the upper
left $2\times 7$ block of the matrix~\eqref{Eq:3x9}.
The $A$-discriminant (Remark \ref{rem:gkz}) for $g(z,q)$ is homogeneous of degree twelve with 127 terms,
\begin{eqnarray*}
  \Delta_{A} &=&q_0^6 q_1 q_2 q_3 q_4 q_5 q_6 - q_0^5 q_1^2 q_2 q_3 q_4^2 q_6 - q_0^5 q_1^2 q_2 q_4 q_5 q_6^2
     - q_0^5 q_1 q_2^2 q_3^2 q_4 q_5 - q_0^5 q_1 q_2^2 q_3 q_5^2 q_6 \\
   &&- q_0^5 q_1 q_3 q_4^2 q_5 q_6^2 - q_0^5q_2  q_3^2 q_4 q_5^2 q_6 + q_0^4 q_1^3 q_2 q_4^2 q_6^2  
     + q_0^4 q_1^2 q_2^2 q_3^2 q_4^2  + \cdots + 16q_3^2 q_4^2 q_5^4 q_6^4\,.
\end{eqnarray*}
This polynomial vanishes at those points $q \in \PP^6$ for which $g(z,q)$ defines a singular curve in $(\CC^\times)^2$.
This six-dimensional family of Laurent polynomials contains a three-dimensional reciprocal subfamily \defcolor{$L_\rec$} obtained by setting
  \begin{equation} \label{eq:substoreciprocal}
    q_0 = p_0, \quad  2q_1 = 2q_2 = p_1, \quad  2q_3 = 2q_4 = p_2, \quad 2q_5 = 2q_6 = p_3\,. 
 \end{equation}
We obtain the following instance of the general reciprocal polynomial~\eqref{eq:fzp}.
 \begin{equation} \label{eq:fzphex} 
   f(z,p) \ = \
   p_0 + p_1 \left( \frac{x + x^{-1}}{2} \right) + p_2 \left( \frac{y + y^{-1}}{2} \right) + p_3 \left( \frac{xy^{-1} + x^{-1}y}{2} \right).
 \end{equation}
Substituting \eqref{eq:substoreciprocal} in $\Delta_{A}$ we obtain 
\[
  \frac{1}{256}  (p_0 + p_1 - p_2 - p_3)  (p_0 - p_1 + p_2 - p_3)  (p_0 - p_1 - p_2 + p_3) (p_0 + p_1 + p_2 + p_3) \cdot (\Delta_A^\rec)^2\,,
\]
where $\Delta_A^\rec = 2 p_0 p_1 p_2 p_3 - p_1^2 p_2^2 - p_1^2 p_3^2 - p_2^2 p_3^2$ is the quartic polynomial
in~\eqref{Eq:rec_A-discrim}.
By definition, the primitive square-free part of this expression is the discriminant $\Delta_r(A)$.
The linear factors are predicted by Lemma \ref{lem:linfacdeltarec}.

The quartic  reciprocal $A$-discriminant $\Delta_A^\rec$ is 
the defining equation of the dual surface of the unirational Chebyshev surface in $\PP^3$ parameterized by 
 \[
   (x,y)\ \longmapsto\ \left( 1: \frac{x + x^{-1}}{2}:\frac{y + y^{-1}}{2} :\frac{xy^{-1} + x^{-1}y}{2} \right)\,.
 \]
 The dual surface~$\Delta_A^\rec = 0$ is parameterized by $(m_{234}: -m_{134}: m_{124}: -m_{123})$, where $m_{ijk}$ are the maximal
 minors of the following $3 \times 4$-matrix: 
 \[
     \begin{pmatrix}
              2 &  x+x^{-1}  & y+y^{-1} &  xy^{-1}+yx^{-1} \\ 
              0 & 1-x^{-2}  &     0    &  y^{-1}-yx^{-2} \\ 
              0 &     0    & 1-y^{-2}  &  x^{-1}-xy^{-2}
     \end{pmatrix}\ .
\]
Further specializing the $A$-discriminant to the dispersion polynomial $D_p$ is done by setting 
\begin{equation} \label{eq:psubs_hexagon}
    p_0\ =\  (\lambda- u)(\lambda- v)-a^2-b^2-c^2\,, \quad p_1\ =\ -2ab\,, \quad p_2\ =\ -2ac\,, \quad p_3\ =\ -2bc\,.
\end{equation}
These are the coefficients of $D_p$ when represented in the form \eqref{eq:fzp}.
With this substitution, the reciprocal discriminant $\Delta_A^\rec$ specializes to $-16\,  a^2b^2c^2 \,  (\lambda- u)  (\lambda- v)$. 
\end{Example}
%%%%%%%%%%%%%%%%%%%%%%%%%%%%%%%%%%%%%%%%%%%%%%%%%%%%%%%%%%%%%%%%%%%%%%%%%%%%%%%%%%%%%%%%%%%%%%%%%%%%

We study this discriminant further in Section \ref{S:Hexagonal}.
As observed in Example \ref{ex:reciprocal_disc_hexagonal}, when $A$ is symmetric in that $A=-A$,
the reciprocal $A$-discriminant is a factor of the specialization of the $A$-discriminant.
This is well-known for univariate reciprocal polynomials, when $d=1$.
Such a polynomial $f$ has the form $f(z)=g(\frac{z+z^{-1}}{2})$, for $g$ a univariate polynomial.
Then
\[
    \disc(f)\ =\ 2^{-2n(n-1)}g(1)\cdot g(-1)\cdot \bigl(\disc(g)\bigr)^2\,,
\]
where $n$ is the degree of $g$.
(The exponent 2 may be deduced from the formula for the discriminant in terms of its roots.)
This suggests that the reciprocal $A$-discriminant is much simpler than the $A$-discriminant, with a significantly smaller degree.
Polynomials with support $A$ form a projective space $\PP^{2m}$, where $A$ consists of $2m+1$ exponent vectors.
The set of reciprocal polynomials forms an $m$-dimensional linear subspace $\defcolor{L_\rec}$. 

%%%%%%%%%%%%%%%%%%%%%%%%%%%%%%%%%%%%%%%%%%%%%%%%%%%%%%%%%%%%%%%%%%%%%%%%%%%%%%%%%%%%%%%%%%%%%%%%%%%%
\begin{theorem}
If $a_1, \ldots, a_m$ generate $\mathbb{Z}^d$ as a group, then $\deg(\Delta_{A}^\rec) \leq \frac{1}{2} (\deg(\Delta_{A})- 2^d)$.
\end{theorem}
%%%%%%%%%%%%%%%%%%%%%%%%%%%%%%%%%%%%%%%%%%%%%%%%%%%%%%%%%%%%%%%%%%%%%%%%%%%%%%%%%%%%%%%%%%%%%%%%%%%%
\begin{proof}
  Recall that the discriminant variety $\nabla_r(A)\subset L_\rec$ for the family~\eqref{eq:fzp} is a hypersurface.
  As $\nabla_r(A)=L_\rec\cap\nabla_A$, we have that $L_\rec\not\subset \nabla_A$.
  Consequently, $\nabla_r(A)$ is a variety (reduced scheme) defined by the restriction
 \defcolor{$\Delta_{A}(L_\rec)$} of the $A$-discriminant to $L_\rec$ and 
 the polynomial $\Delta_r(A)$ defining  $\nabla_r(A)$ is the primitive squarefree part of $\Delta_{A}(L_\rec)$.
 
The assumption that $a_1, \ldots, a_m$ generate $\mathbb{Z}^d$ as a group implies that the hyperplanes/linear factors from Lemma \ref{lem:linfacdeltarec} are distinct, so that $\Delta_{A}(L_\rec)$ has $2^d$ linear factors.
 By construction, one of the remaining factors of $\Delta_A(L_\rec)$ is $\Delta_A^\rec$.
 (This includes the case that $\Delta_A^\rec=1$.)
 We claim that the factor $\Delta_A^\rec$ appears with an exponent $\geq 2$ (as we saw in Example~\ref{ex:reciprocal_disc_hexagonal}).
 Indeed, for a generic point on $\nabla_A^\rec$, the hypersurface \eqref{eq:fzp} has at least two singular points related by taking
 reciprocals.
 Thus the $A$-discriminant has a self-intersection along the reciprocal $A$-discriminant.
\end{proof}
%%%%%%%%%%%%%%%%%%%%%%%%%%%%%%%%%%%%%%%%%%%%%%%%%%%%%%%%%%%%%%%%%%%%%%%%%%%%%%%%%%%%%%%%%%%%%%%%%%%%

We illustrate this theorem in the next example. 

%%%%%%%%%%%%%%%%%%%%%%%%%%%%%%%%%%%%%%%%%%%%%%%%%%%%%%%%%%%%%%%%%%%%%%%%%%%%%%%%%%%%%%%%%%%%%%%%%%%%
\begin{Example}\label{Ex:medDiscriminant}
We add four more columns to the matrix $A$ from Example \ref{ex:reciprocal_disc_hexagonal}.
 \[
    A \ = \ \left(\begin{array}{rrrrrrrrrrr}
0 & 1 & -1 & 0 &  0 &  1 & -1 & 1 & -1 & 2 & -2 \\ 
0 & 0 &  0 & 1 & -1 & -1 &  1 & 1 & -1 & 0 &  0  
  \end{array}\right)\,.
 \]
 Using the software from \cite{rose2023tropical}, we find that the $A$-discriminant has degree 24 and its Newton polytope
 contains 200249 lattice points.
 Computing its coefficients is infeasible, but we may compute the reciprocal discriminant without first computing the $A$-discriminant.
 Let 
 \[ 
  f(z,p) \ = \ 
  p_0 + p_1 \left( \frac{x + x^{-1}}{2} \right) + \cdots  + p_4 \left( \frac{xy + x^{-1}y^{-1}}{2} \right)
    + p_5  \left( \frac{x^2+ x^{-2}}{2} \right)
 \]
 be the polynomial in \eqref{eq:fzphex} plus two new terms.
 Specializing the $A$-discriminant to this reciprocal family, we expect four linear factors (Lemma \ref{lem:linfacdeltarec}) multiplied by
 $(\Delta_A^\rec)^2$.
 From this observation, we expect that $\Delta_A^\rec$ has degree 10.
 We use Algorithm~\ref{alg:discr} to compute $\Delta_r(A)$, taking less than ten seconds and stopping at Step 4.
 It has the four expected linear factors and one of degree 10 with (only) 249 terms, which is $\Delta_A^\rec$.
 We will relate this to a graph in Section \ref{S:four}. 
\end{Example}
%%%%%%%%%%%%%%%%%%%%%%%%%%%%%%%%%%%%%%%%%%%%%%%%%%%%%%%%%%%%%%%%%%%%%%%%%%%%%%%%%%%%%%%%%%%%%%%%%%%%

%%%%%%%%%%%%%%%%%%%%%%%%%%%%%%%%%%%%%%%%%%%%%%%%%%%%%%%%%%%%%%%%%%%%%%%%%%%%%%%%%%%%%%%%%%%%%%%%%%%%
\begin{Example}\label{Ex:bigDiscriminant}
  We add the two columns $(\begin{smallmatrix}0&0\\2&-2\end{smallmatrix})$ to the matrix $A$ of Example~\ref{Ex:medDiscriminant}.
  Let 
 \[ 
  f(z,p) \ = \ 
  p_0 + p_1 \left( \frac{x + x^{-1}}{2} \right) + \cdots  
    + p_5  \left( \frac{x^2+ x^{-2}}{2} \right)  + p_6  \left( \frac{y^2+ y^{-2}}{2} \right)\,.
 \]
  We  use Algorithm~\ref{alg:discr} to compute $\Delta_r(A)$, which takes one hour to compute with Singular and again
  stops at Step 4.
  It has the four expected linear factors and a fifth factor of degree 14 with 3755 terms, 
 \[
   16384 p_0^6 p_3^2 p_4^2 p_5^2 p_6^2\ -\ 65536 p_0^6 p_3^2 p_5^3 p_6^3\ +\ \dotsb\ -\ 524288 p_5^5 p_6^9\ +\ 262144 p_5^4 p_6^{10}\,.
 \]
 We emphasize that this is the unique polynomial in $p_0,\dotsc,p_6$ in a degree reduced Gr\"obner basis for the ideal generated by
 the numerators of  $f, \frac{\partial}{\partial x}f, \frac{\partial}{\partial y}f, xz-1, yw-1$.
 It was not computed using standard techniques for elimination.
\end{Example}
%%%%%%%%%%%%%%%%%%%%%%%%%%%%%%%%%%%%%%%%%%%%%%%%%%%%%%%%%%%%%%%%%%%%%%%%%%%%%%%%%%%%%%%%%%%%%%%%%%%%

%%%%%%%%%%%%%%%%%%%%%%%%%%%%%%%%%%%%%%%%%%%%%%%%%%%%%%%%%%%%%%%%%%%%%%%%%%%%%%%%%%%%%%%%%%%%%%%%%%%%
\section{Resultants and Hessians for bipartite graphs}\label{Sec:bipartite}

Dispersion polynomials are reciprocal in their $z$-variables.
By Lemma~\ref{L:singularFermi}, specializing the discriminant $\Delta_r(A)$ to the coefficients of a dispersion polynomial describes
the dependence of the critical energies on the parameters $p = (e,V)$. 
For bipartite graphs with $|W|=2$, the dispersion polynomial is a determinant of a structured $2 \times 2$ matrix.

Let $H$ be a $2 \times 2$-matrix of the  form, 
 \begin{equation} \label{eq:Lgeneral}
   H \ = \ \begin{pmatrix}
             g(p) & f(z,p) \\ f(z^{-1},p) & h(p)
           \end{pmatrix}, 
 \end{equation}
where $g,h$ are polynomials in $p = (p_1, \ldots, p_n)$,
and $f$ is a bivariate Laurent polynomial in $z$ whose exponent vectors span $\RR^2$ and whose coefficients are polynomials in $p$. We
assume that $f(z,p)=f(z^{-1},p)=0$ has a positive finite number of solutions in $(\CC^\times)^2$, for general $p$.
Let
\begin{equation} \label{eq:D}
    D\ \vcentcolon=\ \det(\lambda {\rm id}_2-H)\ =\
      (\lambda - g(p))(\lambda - h(p)) -  f(z,p) f(z^{-1},p)
\end{equation}
be the characteristic polynomial of $H$.
We consider the incidence variety,   
 \begin{equation}  \label{Eq:Bipartite_Incidence}
  Y\ \vcentcolon=\  \Bigl\{(\lambda, p, z) \in \CC \times \CC^n \times (\CC^*)^2
   \ \big|\  D \ =\ z_1\frac{\partial D}{\partial z_1} \ = \ z_2\frac{\partial D}{\partial z_2} \  = \ 0   \Bigr\} \,.
 \end{equation}
 %

%%%%%%%%%%%%%%%%%%%%%%%%%%%%%%%%%%%%%%%%%%%%%%%%%%%%%%%%%%%%%%%%%%%%%%%%%%%%%%%%%%%%%%%%%%%%%%%%%%%%
\begin{proposition} \label{prop:obviousfactors}
 The closure  $\nabla$ of the projection of $Y$~\eqref{Eq:Bipartite_Incidence} to $\CC \times \CC^n$ contains the hypersurface
 \[
 \Bigl\{(\lambda, p) \;\bigl|\; (\lambda - g(p))(\lambda - h(p))
    \prod_{z^* \in \CP} \bigl((\lambda - g(p))(\lambda - h(p)) - f(z^*,p)^2\bigr)\   =\  0\Bigr\}\,,
 \]
 where $\CP = \{(\pm 1, \pm 1)\}$ is the set of corner points in $(\CC^*)^2$.
\end{proposition}
%%%%%%%%%%%%%%%%%%%%%%%%%%%%%%%%%%%%%%%%%%%%%%%%%%%%%%%%%%%%%%%%%%%%%%%%%%%%%%%%%%%%%%%%%%%%%%%%%%%%

This variety $\nabla$ is a discriminant variety for the family of polynomials \eqref{eq:D}.
Writing $A$ for the exponents of $z$ in $D$, the specialization of the discriminant $\Delta_r(A)$ at the coefficients of
$z$ in $D$ vanishes on $\nabla$.
Proposition \ref{prop:obviousfactors} gives six components/factors of this specialization with ten critical
energies.
The factors which are quadratic in $\lambda$ correspond to the linear factors of  $\Delta_r(A)$ from Lemma~\ref{lem:linfacdeltarec}.
We will see in Section~\ref{S:four} that there may be additional~factors.

%%%%%%%%%%%%%%%%%%%%%%%%%%%%%%%%%%%%%%%%%%%%%%%%%%%%%%%%%%%%%%%%%%%%%%%%%%%%%%%%%%%%%%%%%%%%%%%%%%%%
\begin{proof}[Proof of Proposition \ref{prop:obviousfactors}]
  When $g(p)= \lambda$, we have $D=\det H = -f(z,p)f(z^{-1},p)$.
  The equations $D = z_1\partial_{z_1} D = z_2\partial_{z_2} D = 0$ hold at the singular points of the curve
  defined by $f(z,p)f(z^{-1},p)=0$, which include points where the two branches meet, e.g.\ points $z$ satisfying
  $f(z,p) = f(z^{-1},p) = 0$.
  By assumption, this has solutions and such singular points exist.
  The same argument applies for $h(p)= \lambda$.
  The quadratic factors $D(p, \lambda,z^*)$ were identified in Lemma \ref{lem:linfacdeltarec}.
\end{proof}
%%%%%%%%%%%%%%%%%%%%%%%%%%%%%%%%%%%%%%%%%%%%%%%%%%%%%%%%%%%%%%%%%%%%%%%%%%%%%%%%%%%%%%%%%%%%%%%%%%%%

%%%%%%%%%%%%%%%%%%%%%%%%%%%%%%%%%%%%%%%%%%%%%%%%%%%%%%%%%%%%%%%%%%%%%%%%%%%%%%%%%%%%%%%%%%%%%%%%%%%%
\begin{Example}
  The Floquet matrix $H_p(x,y)$ from Example \ref{ex:floquethexagonal} has the form \eqref{eq:Lgeneral} with
  $f(z,p) = a + bx + cy$, $g(p) =  u$,  and $h(p) =  v$.
  The factors $\lambda -  u$ and $\lambda -  v$ were identified as factors of the specialized reciprocal
  $A$-discriminant in Example \ref{ex:reciprocal_disc_hexagonal}.
\end{Example}
%%%%%%%%%%%%%%%%%%%%%%%%%%%%%%%%%%%%%%%%%%%%%%%%%%%%%%%%%%%%%%%%%%%%%%%%%%%%%%%%%%%%%%%%%%%%%%%%%%%%

We view the defining equation of $\nabla$ as a univariate polynomial in $\lambda$. Its roots are the critical energies. We are interested in parameters $p$ for which two critical energies collide. To this end, for $z^*\in\CP$, we solve $D(\lambda,z^*,p)=(\lambda - g(p))(\lambda - h(p)) - f(z^*,p)^2=0$ for $\lambda$: 
 \begin{equation}\label{Eq:lambda_corner_solve}
   \lambda\ =\ \frac{1}{2}\left( g(p)+h(p)\  \pm \sqrt{(g(p)-h(p))^2 +  4\,f(z^*,p)^2}\right)\ .
 \end{equation}
Observe that  $g(p)$ is a root of $D(p,\lambda, z^*)$ if and only if $f(z^*,p)=0$.

%%%%%%%%%%%%%%%%%%%%%%%%%%%%%%%%%%%%%%%%%%%%%%%%%%%%%%%%%%%%%%%%%%%%%%%%%%%%%%%%%%%%%%%%%%%%%%%%%%%%
 \begin{proposition}
   \label{prop:Disc_Resultant}
  Let $D, f,g,h$ be as in Proposition \ref{prop:obviousfactors}, and $z^*,z_1^*, z_2^* \in \CP$.
  \begin{enumerate}
  \item The discriminant of $D(\lambda,z^*,p)$ is
    $\nabla_\lambda (D(\lambda,z^*,p)) \ = \ (g(p)-h(p))^2 + 4  f(z^*,p)^2$.\vspace{1pt}
   \item The resultant ${\rm res}_{\lambda}(D(p,\lambda,z_1^*), D(p,\lambda,z_2^*))$ is $(f(z_1^*,p)^2 - f(z_2^*,p)^2)^2$.\vspace{1pt}
   \item The resultant ${\rm res}_{\lambda}(D(p,\lambda,z^*), g(p)-\lambda)$ is $-f(z^*,p)^2$.
\end{enumerate}
\end{proposition}
%%%%%%%%%%%%%%%%%%%%%%%%%%%%%%%%%%%%%%%%%%%%%%%%%%%%%%%%%%%%%%%%%%%%%%%%%%%%%%%%%%%%%%%%%%%%%%%%%%%%
\begin{proof}
  The first item is part of the calculation of~\eqref{Eq:lambda_corner_solve}, and the others follow from
  the Poisson formula for the resultant of univariate monic  polynomials $F(t),G(t)$, 
  \[
    {\rm res}(F(t),G(t))\ =\ \pm \prod_{a\colon F(a)=0}G(a)\,.  \qedhere
  \]
\end{proof}
%%%%%%%%%%%%%%%%%%%%%%%%%%%%%%%%%%%%%%%%%%%%%%%%%%%%%%%%%%%%%%%%%%%%%%%%%%%%%%%%%%%%%%%%%%%%%%%%%%%%

We derive a formula for the Hessian matrix of the eigenvalue function $\lambda(z) = \lambda(x,y)$ evaluated at the corner points $z \in \CP$.
This will be convenient to describe the nature of the corresponding critical points in later sections.
We suppress the dependence on $p$ to simplify notation and start from the implicit definition of $\lambda(x,y)$,
\[
     D(x,y,\lambda(x,y)) \  = \  (\lambda(x,y)-g)(\lambda(x,y)-h) - f(x,y)f(x^{-1},y^{-1}) \ = \ 0\,.
\]
Writing $\lambda_x = \frac{\partial \lambda(x,y)}{\partial x}$ and similarly for $y$, implicit differentiation of $D=0$ gives
  \begin{equation} \label{eq:dlambdadx}
    2 \lambda \lambda_x - (g+h) \lambda_x \ = \ f_x(x,y)f(x^{-1},y^{-1}) - f(x,y)f_x(x^{-1},y^{-1}) \frac{1}{x^2}
  \end{equation}
Taking $\frac{\partial }{\partial x}$ once more and evaluating at a point $(x,y) \in \CP$, for which $x^2 = y^2 = 1$ and
$\lambda_x(x,y) = \lambda_y(x,y) = 0$, we obtain 
 \[
   \Big (\lambda - \frac{g+h}{2} \Big ) \lambda_{xx} \  = \   f(x,y) f_{xx}(x,y) - f_x(x,y)^2 + xf(x,y)f_x(x,y)\,.
 \]
A similar calculation for  $\frac{\partial }{\partial y}$ gives
\[
   \Big (\lambda - \frac{g+h}{2} \Big ) \lambda_{xy} \  = \  f(x,y) f_{xy}(x,y) - f_x(x,y) f_y(x,y)\,.
\] 
This gives an expression for the Hessian matrix of $\lambda(x,y)$ evaluated at $(x,y) \in \CP$,
\begin{equation} \label{eq:hess}
  \Big( \lambda - \frac{g+h}{2} \Big ){\rm Hess}(\lambda) \  = \
    \begin{pmatrix}
      f f_{xx} - f_x^2 + xf f_x & f f_{xy}-f_xf_y \\ 
         f f_{xy}-f_xf_y       & f f_{yy} -f_y^2 + yf f_y
    \end{pmatrix}\ ,
\end{equation}
where $f, f_x, f_y, \ldots$ are the values of $f$ and its partial derivatives at $(x,y)$.
We will use~\eqref{eq:hess} to investigate the nature of the critical points $(x,y) \in \CP$.
The sign of the factor $(\lambda - \tfrac{1}{2}(g+h))$ is determined by which  root~\eqref{Eq:lambda_corner_solve}
of $D$ we are considering.

%%%%%%%%%%%%%%%%%%%%%%%%%%%%%%%%%%%%%%%%%%%%%%%%%%%%%%%%%%%%%%%%%%%%%%%%%%%%%%%%%%%%%%%%%%%%%%%%%%%%
\begin{Example}\label{Ex:Hex_Hess}
If $f = a + bx + cy$, then \eqref{eq:hess} becomes
\[
  \begin{pmatrix}
     -b^2 + x(a + bx + cy)b & -bc \\ 
    -bc & -c^2 + y(a + bx + cy)c
  \end{pmatrix} \  = \
  \begin{pmatrix}
   abx +bcxy & -bc \\ -bc & acy + bcxy
  \end{pmatrix}\ ,
\]
using that $x^2 = y^2 = 1$ for the simplification. 
\end{Example}
%%%%%%%%%%%%%%%%%%%%%%%%%%%%%%%%%%%%%%%%%%%%%%%%%%%%%%%%%%%%%%%%%%%%%%%%%%%%%%%%%%%%%%%%%%%%%%%%%%%%

%%%%%%%%%%%%%%%%%%%%%%%%%%%%%%%%%%%%%%%%%%%%%%%%%%%%%%%%%%%%%%%%%%%%%%%%%%%%%%%%%%%%%%%%%%%%%%%%%%%%
\section{The Hexagonal Lattice}\label{S:Hexagonal}
Consider the family of operators on the hexagonal lattice of Figure~\ref{fig:HexGraph}.
Its Floquet matrix $H_p(x,y)$ was given in Example~\ref{ex:floquethexagonal} and its dispersion polynomial $D_p(x,y, \lambda)$ is
 \begin{equation}\label{Eq:Hex_dispersion}
  (\lambda -  u)(\lambda -  v)\ -\ (a+bx+cy)(a+bx^{-1}+cy^{-1})\,.
 \end{equation}
Eliminating the variables $x$ and $y$ from the critical point equations~\eqref{eq:critpoints}
gives a degree ten polynomial in $\lambda$.
It has the six irreducible factors as predicted by Proposition \ref{prop:obviousfactors}.
\begin{equation} \label{eq:6factors} \begin{matrix}
  \lambda -  u \\
  \lambda -  v \\
  (\lambda -  u)(\lambda -  v) - (a + b + c)^2 \\
  (\lambda -  u)(\lambda -  v) - (a + b - c)^2 \\
  (\lambda -  u)(\lambda -  v) - (a - b + c)^2 \\
  (\lambda -  u)(\lambda -  v) - (a - b - c)^2 
\end{matrix}
\end{equation}
By Lemma \ref{L:singularFermi}, these give the relations between the parameters $p$ and the critical energies.
For generic values of $p = (a,b,c, u, v)$, there are ten distinct (complex) critical energies corresponding to
the ten roots of the polynomials in \eqref{eq:6factors}.
This fails when a discriminant or a resultant from Proposition \ref{prop:Disc_Resultant} vanishes, which is a necessary condition
for two critical points $(z,\lambda)$ to collide.
Hence, these conditions yield a Bloch discriminant.  \medskip

%%%%%%%%%%%%%%%%%%%%%%%%%%%%%%%%%%%%%%%%%%%%%%%%%%%%%%%%%%%%%%%%%%%%%%%%%%%%%%%%%%%%%%%%%%%%%%%%%%%%
\noindent {\bf Theorem~\ref{T:hexagonal}.}  {\it
  Let $\Gamma$ be the hexagonal graph from Figure \ref{fig:HexGraph}.
  The quintic polynomial
 \[
    \BD_\Gamma \ \vcentcolon=\  (u-v)(a+b+c)(a+b-c)(a-b+c)(a-b-c)\,.
 \]
 is a Bloch discriminant for $\Gamma$. }\smallskip
%%%%%%%%%%%%%%%%%%%%%%%%%%%%%%%%%%%%%%%%%%%%%%%%%%%%%%%%%%%%%%%%%%%%%%%%%%%%%%%%%%%%%%%%%%%%%%%%%%%%

The vanishing of any two of the last four factors implies that one of $a,b,c$ also vanishes, which is forbidden.

%%%%%%%%%%%%%%%%%%%%%%%%%%%%%%%%%%%%%%%%%%%%%%%%%%%%%%%%%%%%%%%%%%%%%%%%%%%%%%%%%%%%%%%%%%%%%%%%%%%%
\begin{proof}
  The resultants from Proposition \ref{prop:Disc_Resultant}(2) do not contribute to the Bloch discriminant, 
  as the corresponding critical points lie above different corner points.
  The resultants from Proposition \ref{prop:Disc_Resultant}(3) are (the squares of) $(a \pm b \pm c)$.

The discriminants from Proposition \ref{prop:Disc_Resultant}(1) are the sums of squares 
\[
    ( u -  v)^2 + 4 \, (a \pm b \pm c)^2 \, .
\]
These vanish only when $ u -  v = a \pm b \pm c = 0$, which is accounted for by $\BD_\Gamma$.
\end{proof}
%%%%%%%%%%%%%%%%%%%%%%%%%%%%%%%%%%%%%%%%%%%%%%%%%%%%%%%%%%%%%%%%%%%%%%%%%%%%%%%%%%%%%%%%%%%%%%%%%%%%

We investigate how the nature of the critical points depends on the parameters.
The eliminant in \eqref{eq:6factors} splitting into many irreducible factors of low degree facilitates this detailed analysis.
The roots of each factor, as a polynomial in $\lambda$, are the critical energies.
We determine the $(x,y)$-coordinates of the corresponding critical point(s), investigate their reality, and the signature of the Hessian
matrix of the energy function $\lambda(x,y)$.
We treat the four quadratic and the two linear factors in \eqref{eq:6factors}~separately.
Sections \ref{subsec:hexagonquadratic} and \ref{subsec:hexagonlinear} imply the following result.

%%%%%%%%%%%%%%%%%%%%%%%%%%%%%%%%%%%%%%%%%%%%%%%%%%%%%%%%%%%%%%%%%%%%%%%%%%%%%%%%%%%%%%%%%%%%%%%%%%%%
\begin{theorem}\label{T:Critical_points_hex}
  If $\BD_\Gamma\neq 0$, then the complex Bloch variety has twelve nondegenerate critical points.
  If the squares $(a \pm b \pm c)^2$ are distinct, then the critical points are distributed among the roots of the factors
  in~\eqref{eq:6factors} as follows. 
  \begin{enumerate}
    \item Each critical energy from a quadratic factor corresponds to one critical point.
    \item Each critical energy from a linear factor corresponds to two critical points.
  \end{enumerate}
\end{theorem}

We note that, if the squares $(a \pm b \pm c)^2$ are not distinct, then some of the quadratic factors in
\eqref{eq:6factors} coincide.
Their critical energies $\lambda$ coincide as well, but the $(x,y)$-coordinates of the corresponding critical points are
distinct: they are distinct corner points.
Therefore, the condition that the squares $(a \pm b \pm c)^2$ are distinct does not appear in the Bloch discriminant.
%%%%%%%%%%%%%%%%%%%%%%%%%%%%%%%%%%%%%%%%%%%%%%%%%%%%%%%%%%%%%%%%%%%%%%%%%%%%%%%%%%%%%%%%%%%%%%%%%%%%

%%%%%%%%%%%%%%%%%%%%%%%%%%%%%%%%%%%%%%%%%%%%%%%%%%%%%%%%%%%%%%%%%%%%%%%%%%%%%%%%%%%%%%%%%%%%%%%%%%%%
\subsection{Critical points for quadratic factors} \label{subsec:hexagonquadratic}

The four quadratic factors~\eqref{eq:6factors} have roots 
 \begin{equation}\label{Eq:Hexagonal_energies}
   \frac{1}{2} \Bigl(u+v\  \pm \sqrt{(u-v)^2 + 4(a \pm b \pm c)^2 }\,\Bigr)\,,
 \end{equation}
see \eqref{Eq:lambda_corner_solve}.
When the squares $(a\pm b \pm c)^2$ are distinct, they give eight distinct critical energies.

For each such critical energy $\lambda$, the Fermi curve $\calF(\lambda)$ is singular at the critical point
$(\pm1,\pm1,\lambda)$ whose signs correspond to the signs in the term $(a \pm b \pm c)$ in~\eqref{Eq:Hexagonal_energies}.
The Newton polytope of the Fermi curve is the hexagon which is the base of the pyramid in  Figure~\ref{Fig:hexPolytope}.
It has one integer point in the interior, so $\calF(\lambda)$ has arithmetic genus
one and is a rational curve with a unique singular point.

We determine the type of critical points above corner points $(x^*,y^*)=(\pm1,\pm1)$.
The Bloch variety has two branches $\lambda^+$ and $\lambda^-$ (the signs corresponding to the sign in front of the square root
in Equation \eqref{Eq:Hexagonal_energies})---these are the spectral band functions.
A point $(x,y,\lambda^+)$ is a relative maximum/saddle point/minimum if and only if $(x,y,\lambda^-)$ is a
relative minimum/saddle point/maximum.
The next proposition uses the notation $\Delta = {\rm BD}_{\Gamma}/(u-v)$. 

\begin{figure}[htb]
\centering
\raisebox{-37pt}{
  \begin{picture}(87,80)(-8,-10)
    \put(0,0){\includegraphics{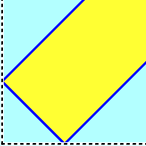}}
    %\put(-8,25){{\color{red}Image not pushed}}
    \put(40,37){I}   \put(53 ,6){II}    \put(5,55){III}   \put(3,5){IV}     
    \put(-8,62){$\frac{c}{a}$}  \put(73,-2){$\frac{b}{a}$}
    \put(-7,28){$1$}     \put(29,-10){$1$} 
  \end{picture}}
\qquad\quad
\begin{tabular}{|r|c|c|c|c|}\hline
  &(1,1)&(-1,1)&(1,-1)&(-1,-1)\\\hline
  \mbox{I}  & Maximum & Saddle  & Saddle  & Saddle  \\\hline
  \mbox{II} & Maximum & Minimum & Saddle  & Saddle  \\\hline
  \mbox{III}& Maximum & Saddle  & Minimum & Saddle  \\\hline
  \mbox{IV} & Maximum & Saddle  & Saddle  & Minimum \\\hline
\end{tabular}
\caption{Critical point analysis for $\lambda^+$ and $a,b,c>0$.\rule{0pt}{14pt} }\label{F:CPA}
\end{figure}

%%%%%%%%%%%%%%%%%%%%%%%%%%%%%%%%%%%%%%%%%%%%%%%%%%%%%%%%%%%%%%%%%%%%%%%%%%%%%%%%%%%%%%%%%%%%%%%%%%%%
\begin{proposition}\label{P:hexagonal_CPA}
  Suppose that $a,b,c\neq 0$ and $\Delta\neq 0$.
  The branch $\lambda^+$ has a nondegenerate critical point lying above each corner point.
  For $a,b,c>0$, this is summarized in Figure~\ref{F:CPA}.
  For the other quadrants, multiply $(x^*,y^*)$ by the sign of $(\frac{b}{a},\frac{c}{a})$.
  When $\Delta<0$ (region I in Figure~\ref{F:CPA}), the global minimum does not occur at a corner point and
  when $\Delta=0$, there is one critical point of multiplicity three at the
  corresponding corner point. 
\end{proposition}
%%%%%%%%%%%%%%%%%%%%%%%%%%%%%%%%%%%%%%%%%%%%%%%%%%%%%%%%%%%%%%%%%%%%%%%%%%%%%%%%%%%%%%%%%%%%%%%%%%%%
%%%%%%%%%%%%%%%%%%%%%%%%%%%%%%%%%%%%%%%%%%%%%%%%%%%%%%%%%%%%%%%%%%%%%%%%%%%%%%%%%%%%%%%%%%%%%%%%%%%%
\begin{proof}
 Let $f$ be a real bivariate function with critical point $(x^*,y^*)$.
 Let \defcolor{$\frac{1}{2}Qf$} be the quadratic term in the Taylor expansion of $f$ at the point $(x^*,y^*)$.
 Then $(x^*,y^*)$ is a relative extremum if the determinant of $Qf$ is positive and a saddle point if the determinant is negative.
 A relative extremum is a relative minimum, resp.~maximum, when the trace of $Qf$ is positive, resp.~negative.
 Important for us is that $Qf$ is a linear map on the second symmetric power of the tangent space at $(x^*,y^*)$.

 In Calculus, $Qf$ is the matrix $\Hess(f)$ of second partial derivatives of $f$ at $(x^*,y^*)$.
 For us, however,  $Qf=-\Hess(f)$.
 This is because our ambient space is $\TT^2$ and not $\RR^2$; the tangent space $T_x\TT$ to $\TT$ at $x\in\TT$ is
 $\RR\cdot \sqrt{-1}x\subset\CC=T_x\CC^\times$.
 At a  corner point $(x^*,y^*)$, we have
 $T_{(x^*,y^*)}\TT^2=\RR\sqrt{-1} \oplus \RR\sqrt{-1}\subset \CC^2=T_{(x^*,y^*)}(\CC^\times)^2$.
 The $\sqrt{-1}$ gives a negative sign in the second symmetric power.

 Following the discussion before Example~\ref{Ex:Hex_Hess}, as well as~\eqref{eq:hess}, we have 
  \[
  \delta \cdot Q\lambda^+\ =\ -\delta\cdot \Hess(\lambda^+)\ \vcentcolon=\
  \begin{pmatrix} - abx^* - bcx^*y^*  & bc\\ bc & - ac y^* - bcx^*y^* \end{pmatrix}\ ,
  \]
  where  $\delta\vcentcolon=\lambda^+-(u+v)/2 > 0$.
  Thus
  \begin{eqnarray*}
     \delta^2\det(Q\lambda^+)        &=& a bc x^*y^* ( a +bx^* + cy^*)\,, \ \mbox{ and} \\
     \delta \mbox{Trace}(Q\lambda^+) &=& -(abx^* + acy^* + 2bcx^*y^*)\,.
  \end{eqnarray*}
  Write \defcolor{$\det$} for $\delta^2\det(Q\lambda^+)$ and \defcolor{$\Tr$} for $\delta \mbox{Trace}(Q\lambda^+)$.
  We analyze their signs at each corner point $(x^*,y^*)$, assuming that $a,b,c>0$.

  \begin{itemize}
    \item[]\hspace{-38pt}$(+1,+1)$:   Then $\det>0$ and $\Tr<0$.

    \item[]\hspace{-38pt}$(-1,+1)$:   Then $\det=-abc(a-b+c)>0$ if and only if $a-b+c<0$, which is region II in Figure~\ref{F:CPA}.
      Assuming this, as $b-c>a$ we have
      \[
        \Tr\ =\ ab-ac+2bc \ =\ a(b-c)+2bc \ >\ a^2+2bc\ >\ 0\,.
      \]

    \item[]\hspace{-38pt}$(+1,-1)$:  This is symmetric to the previous case, so $\det>0$ if and only if $(a,b,c)$ lies in region III,
      and there $\Tr>0$.

    \item[]\hspace{-38pt}$(-1,-1)$: Then $\det = abc(a-b-c)>0$ if and only if $a>b+c$, which is region IV.
      Assuming that, $\Tr=ab+ac-2bc > (b+c)^2-2bc=b^2+c^2>0$. 
  \end{itemize}

  This verifies the claims in Figure~\ref{F:CPA}.
  The statement about multiplicities is a consequence of our subsequent analysis---in any perturbation about a point
  $a,b,c>0$  with $\Delta=0$, two additional nondegenerate critical points coalesce at the corresponding corner point.
  When $\Delta<0$, we will show there are two minima which are not at a corner point.
  \end{proof}
%%%%%%%%%%%%%%%%%%%%%%%%%%%%%%%%%%%%%%%%%%%%%%%%%%%%%%%%%%%%%%%%%%%%%%%%%%%%%%%%%%%%%%%%%%%%%%%%%%%%

\subsection{Critical points when $\lambda = u$ or $\lambda = v$} \label{subsec:hexagonlinear}
The study of the linear factors differs depending on if $u \neq v$ or $u = v$.
\subsubsection{Analysis of the linear factors in \eqref{eq:6factors} for $u \neq v$.}
If $\lambda =  u$ or $\lambda =  v$, then the determinant $D(z,\lambda)$ factors
and the Fermi variety is reducible,
\begin{equation}\label{Eq:twoFactors}
    D(z,\lambda) \ = \  -(a + bx + cy) (a + bx^{-1} + cy^{-1})\,.
\end{equation}
By Lemma~\ref{L:singularFermi}, the complex critical points at $\lambda$ are singular points of the Fermi
variety $\calF(\lambda)$.
As shown in the proof of Proposition~\ref{prop:obviousfactors}, these are 
the intersections of the curves defined by the two factors~\eqref{Eq:twoFactors}, as each curve is smooth.
For general $a,b,c$, there are two such points.
Indeed, let us solve the system
\[
    a + bx + cy\ =\ a + bx^{-1} + cy^{-1} \ =\ 0\,.
\]
Multiplying by $a+bx^{-1}$ and $cy$ respectively, then subtracting one from the other gives 
 \begin{equation}\label{Eq:UnivConsequence}
   0\ =\ abx + (a^2+b^2-c^2) + abx^{-1}\,.
 \end{equation}
Substituting its two solutions into $y=-\frac{1}{c}(a+bx)$ gives the two critical points
 \begin{equation}\label{Eq:D-points}
   (x,y)\ =\ \Bigl(-\frac{a^2+b^2-c^2}{2ab}  \ \pm\ \frac{1}{2ab}\sqrt{\Delta} \ ,\ 
           -\frac{a^2-b^2+c^2}{2ac}  \ \mp\ \frac{1}{2ac}\sqrt{\Delta}\Bigr)\ ,
 \end{equation}
 where $\Delta$ is as before.
 It is also the discriminant of the polynomial~\eqref{Eq:UnivConsequence},
\[
\Delta\ =\ (a^2+b^2-c^2)^2-4a^2b^2 \ =\ (a - b - c) (c - a - b) (b - a - c) (a + b + c)\,.
\]
The two points~\eqref{Eq:D-points} are interchanged under the involution $(x,y) \mapsto (x^{-1},y^{-1})$.

When $\Delta>0$, the two points~\eqref{Eq:D-points} lie in $\RR^2$ and are therefore not critical points of the real
Bloch variety. 
When $\Delta<0$, they lie on the compact torus $\TT^2$ and both are critical points of the real Bloch variety.
We display this in the $(b,c)$-plane where $a=1$.
\[
  \begin{picture}(125,72)
    \put(0,0){\includegraphics{images/triangleIneq}}
    \put(21,64){\small$c$}  \put(75,15){\small$b$}
    \put(95,35){$\RR^3_+$}
    \put(95,61){$\Delta<0$} \put(93,64){\vector(-1,0){31}}
   \end{picture}
\]
The dotted lines are where $b$ or $c$ vanish, the solid lines are where $\Delta=0$.
The discriminant $\Delta$ is negative in the yellow regions and positive in the light blue regions.

This number, 2, of solutions is the mixed area of the Newton polygons of the factors in~\eqref{Eq:twoFactors}, which is the upper bound
given by Bernstein's Theorem~\cite{Bernstein}.
Consequently, the two curves meet transversally at each solution when $\Delta\neq 0$. 
Below, we show the Newton polytopes $P$ and $-P$ of each factor, as well as their Minkowski sum $P+(-P)$, which
is the hexagonal base of the Newton polytope $\calN$ of the dispersion relation shown in Figure~\ref{Fig:hexPolytope}.
 \begin{equation}\label{Eq:MinkowskiSum}
   P\ =\ \raisebox{-20pt}{\includegraphics{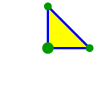}}\qquad
   -P =\ \raisebox{-20pt}{\includegraphics{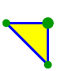}}\qquad
   P+(-P)\ =\ \raisebox{-20pt}{\includegraphics{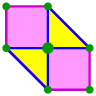}}
 \end{equation}
The mixed area is given by $\mbox{Area}(P+(-P))-\mbox{Area}(P)-\mbox{Area}(-P)=3-\frac{1}{2}-\frac{1}{2}=2$.

As $a,b,c\neq 0$, at most one factor of the discriminant $\Delta$ vanishes. 
This factor is $a\pm b \pm c$ for appropriate signs.
In this case, the two points~\eqref{Eq:D-points} coincide with the corresponding 
corner point $(\pm 1, \pm 1)$ forming a degenerate critical point.
This explains the appearance of $\Delta$ as a factor in the Bloch discriminant
and the multiplicity 3 in Proposition~\ref{P:hexagonal_CPA}.

This also determines the spectrum of $H_p$, showing that it has a band gap.
If $\Delta>0$, so that the critical points from the linear factors of~\eqref{eq:6factors} do not lie in $\TT^2$,
and if we let \defcolor{$m$}, respectively \defcolor{$M$}, be the minimum (respectively maximum) of the quantities
 \[
     \left\{\tfrac{1}{2} \sqrt{(u-v)^2 + 4(a \pm b \pm c)^2 } \mid \mbox{ all choices of $\pm$} \right\}\,,
 \]
then by~\eqref{Eq:Hexagonal_energies}, the spectrum of $H_p$ is
\[
   \left[ \tfrac{u+v}{2}-M, \tfrac{u+v}{2}-m\right]\ \bigcup\ 
   \left[ \tfrac{u+v}{2}+m, \tfrac{u+v}{2}+M\right]\, ,
\]
with a band gap of $2m$.   

If $\Delta<0$, let us assume that $v<u$.
We now have four more critical points, two from each linear factor of~\eqref{eq:6factors}.
As $0<\frac{u-v}{2} < m$,
\[
   \tfrac{u+v}{2}\ <\  u \ =\ \tfrac{u+v}{2} + \tfrac{u-v}{2}\ <\ \tfrac{u+v}{2}+m\,.
\]
We also have $\tfrac{u+v}{2}>v>\tfrac{u+v}{2}-m$.
Thus the spectrum of $H$ is 
\[
   \left[ \tfrac{u+v}{2}-M, v\right]\ \bigcup\ 
   \left[  u, \tfrac{u+v}{2}+M\right]\, ,
\]
which has a band gap of $u-v$.
Thus both critical points at $\lambda=u$ are relative minima and both at $\lambda=v$ are relative maxima.

In either case, the critical point equations have twelve solutions, although not all are real (lie in $\TT^2$).
As this equals Kushnirenko's bound~\cite{Kushnirenko}  (normalized volume of the Newton polytope shown in Figure~\ref{Fig:hexPolytope})
all solutions have multiplicity 1. 
By~\cite[Thm.~5.2]{FS24}, all critical points are nondegenerate and the spectral edges conjecture holds when $u\neq v$ and~$\Delta\neq 0$.

\begin{figure}[htb]
\centering
\begin{picture}(435,143)(-100,-15)

  \put(140,0){\includegraphics{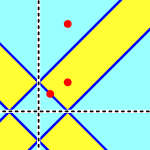}}

  \put(-100,-15){\includegraphics[height=90pt]{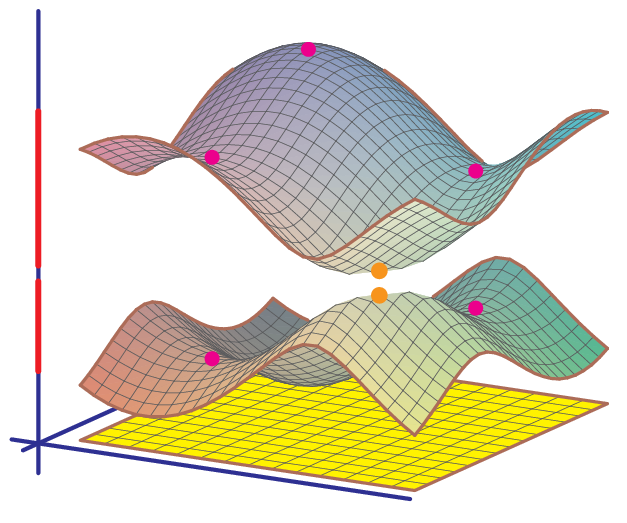}}
  \put( 230,-10){\includegraphics[height=90pt]{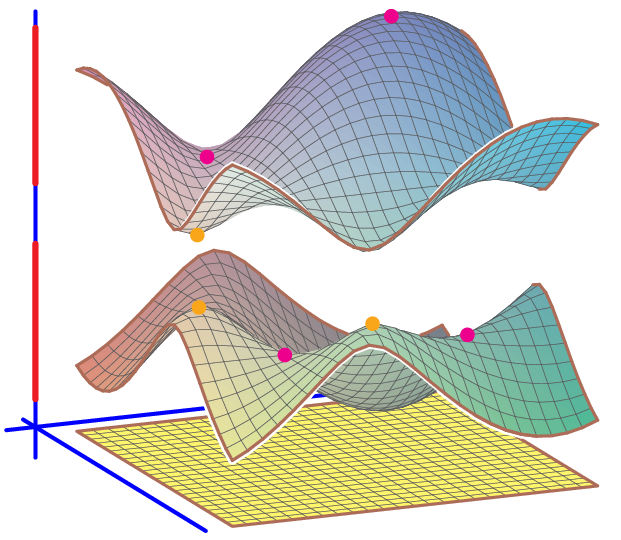}}
  \put(  20, 35){\includegraphics[height=90pt]{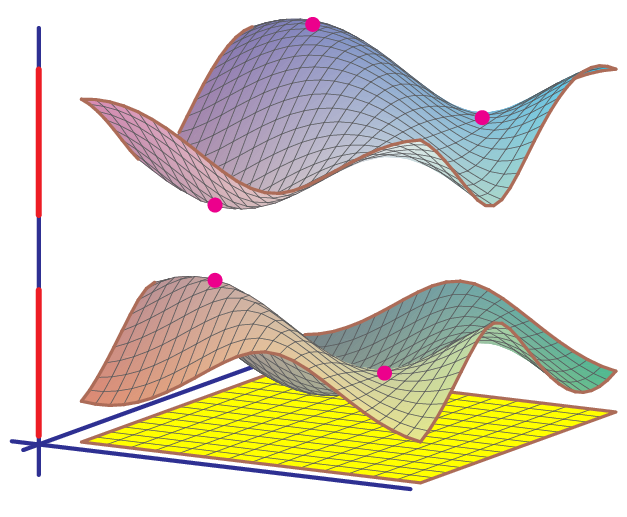}}

  \put(-45, 80){$\Delta=0$}
  \put(133,113){$\Delta>0$}
  \put(265, 85){$\Delta<0$}

  \thicklines
  \put(74,26.5){{\color{white}\vector(1,0){87.5}}}  \put(74,27.5){{\color{white}\vector(1,0){87.5}}}
  \put(74,27){{\color{white}\vector(1,0){87.5}}}
  \put(119.5 ,86.5){{\color{white}\vector(2,-1){50}}}  \put(120,87.5){{\color{white}\vector(2,-1){50}}}
  \put(119.75,87  ){{\color{white}\vector(2,-1){50}}}
  \put(226,32){{\color{white}\vector(-1,0){51}}}  \put(226,33){{\color{white}\vector(-1,0){51}}}
  \put(226,32.5){{\color{white}\vector(-1,0){51}}}
  \thinlines
  \put(14,27){\vector(1,0){146.5}}
  \put(120,87){\vector(2,-1){49}}
  \put(226,32.5){\vector(-1,0){50}}
\end{picture}

\caption{Bloch varieties with $u\neq v$ at parameters $(1,\frac{2}{5},\frac{3}{5})$, $(1,1,3)$, and $(1,1,1)$.}
\label{F:Hex_LandscapeI}
\end{figure}

When $\Delta=0$ and $u\neq v$, exactly one of the factors in $\Delta$ vanishes and as we remarked there are two
degenerate critical points corresponding to that factor at energies $u$ and $v$.
Figure~\ref{F:Hex_LandscapeI} shows three Bloch varieties with their corresponding parameter values.
%%%%%%%%%%%%%%%%%%%%%%%%%%%%%%%%%%%%%%%%%%%%%%%%%%%%%%%%%%%%%%%%%%%%%%%%%%%%%%%%%%%%%%%%%%%%%%%%%%%%
%%%%%%%%%%%%%%%%%%%%%%%%%%%%%%%%%%%%%%%%%%%%%%%%%%%%%%%%%%%%%%%%%%%%%%%%%%%%%%%%%%%%%%%%%%%%%%%%%%%%

%%%%%%%%%%%%%%%%%%%%%%%%%%%%%%%%%%%%%%%%%%%%%%%%%%%%%%%%%%%%%%%%%%%%%%%%%%%%%%%%%%%%%%%%%%%%%%%%%%%%
\begin{table}[htb]
\centering
\begin{tabular}{||c | c | c | l ||} 
 \hline
$\Delta$ &  nondegenerate & degenerate & Hessian \\ [0.3ex] 
 & critical points & critical points &  \\ [0.3ex] 
 \hline\hline 
  $< 0$ & 12 & 0 &\  (3,6,3) \rule{0pt}{12pt}\\  [0.6ex]
  $> 0$ & 8 & 0 &\  (2,4,2) \\ [0.6ex]
  $= 0$ & 6 & 2 (mult. 3) &\ (1,4,1) \\ [0.2ex]
  \hline
\end{tabular}
\caption{Critical point summary\rule{0pt}{14pt} % for the hexagonal lattice for
          when $u \neq v$.}
\label{table:hexagonal_second}
\end{table}
%%%%%%%%%%%%%%%%%%%%%%%%%%%%%%%%%%%%%%%%%%%%%%%%%%%%%%%%%%%%%%%%%%%%%%%%%%%%%%%%%%%%%%%%%%%%%%%%%%%%
     
\subsubsection{Analysis of the linear factors in \eqref{eq:6factors} for $u = v$.}
When $u=v$ the two linear factors in~\eqref{eq:6factors} become equal.
Setting $\kappa:=\lambda-u=\lambda-v$, the dispersion polynomial becomes
 \begin{equation}
   \label{Eq:Hexagon_ODP}
   \kappa^2\ -\  (a + bx + cy)(a + bx^{-1} + cy^{-1})\,.
 \end{equation}
When $\Delta\neq 0$, the curves defined by the two factors meet transversally at the points~\eqref{Eq:D-points}.
Then $\defcolor{w}\vcentcolon= a + bx + cy$ and $\defcolor{z}\vcentcolon= a + bx^{-1} + cy^{-1}$ give local coordinates at these points.
In these coordinates, the Bloch variety is defined by $\kappa^2-wz=0$, showing that these critical points when $u=v$ are
ordinary double point singularities on the complex Bloch variety.

\begin{figure}[htb]
\centering
\begin{picture}(435,138)(-100,-15)

  \put(140,0){\includegraphics{images/triangleIneqDots}}

  \put(-100,-15){\includegraphics[height=90pt]{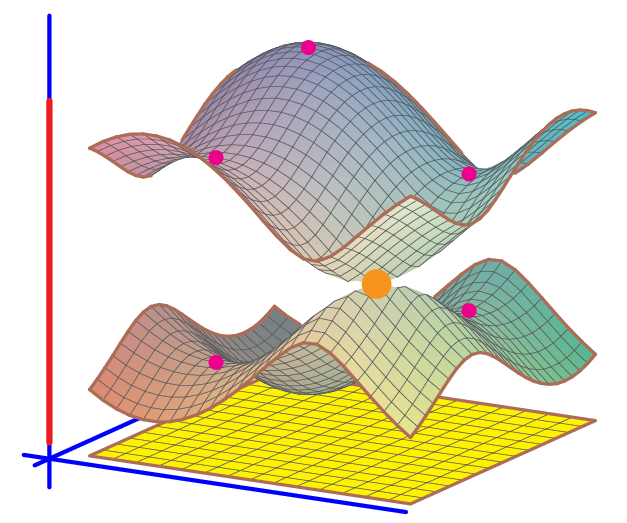}}
  \put( 230,-10){\includegraphics[height=90pt]{images/hexagonal_111}}
  \put(  20, 35){\includegraphics[height=90pt]{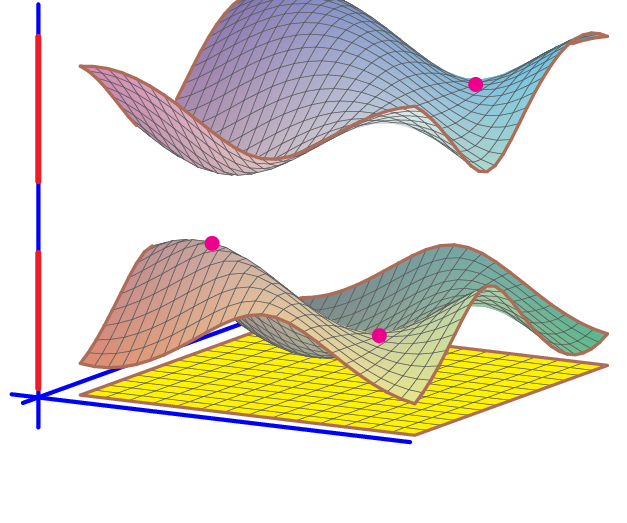}}

  \put(-45, 80){$\Delta=0$}
  \put(138,113){$\Delta>0$}
  \put(265, 85){$\Delta<0$}

  \thicklines
  \put(74,26.5){{\color{white}\vector(1,0){87.5}}}  \put(74,27.5){{\color{white}\vector(1,0){87.5}}}
  \put(74,27){{\color{white}\vector(1,0){87.5}}}
  \put(119.5 ,86.5){{\color{white}\vector(2,-1){50}}}  \put(120,87.5){{\color{white}\vector(2,-1){50}}}
  \put(119.75,87  ){{\color{white}\vector(2,-1){50}}}
  \put(226,32){{\color{white}\vector(-1,0){51}}}  \put(226,33){{\color{white}\vector(-1,0){51}}}
  \put(226,32.5){{\color{white}\vector(-1,0){51}}}
  \thinlines
  \put(14,27){\vector(1,0){146.5}}
  \put(120,87){\vector(2,-1){49}}
  \put(226,32.5){\vector(-1,0){50}}
\end{picture}

\caption{Bloch varieties with $u=v$ at parameters $(1,\frac{2}{5},\frac{3}{5})$, $(1,1,3)$, and $(1,1,1)$.}
\label{F:Hex_LandscapeII}
\end{figure}

If $\Delta>0$, then the  points~\eqref{Eq:D-points} lie in $\RR^2$ and the corresponding critical points are not on the
real Bloch variety. 
Thus, the real Bloch variety is smooth with eight nondegenerate critical points and a band gap as $|a\pm b\pm c|>0$.
We see this at the top (middle) of  Figure~\ref{F:Hex_LandscapeII}.

If $\Delta<0$, then the points~\eqref{Eq:D-points} lie in $\TT^2$.
To analyze the singularity, we introduce local coordinates that are fixed under the non-standard real structure (conjugation) giving
$\TT^2$. 
Observe that under this conjugation, $\overline{w}=z$.
Then $r=\frac{w+z}{2}$ and $s=\frac{w-z}{2\sqrt{-1}}$ are real (fixed under conjugation) local coordinates at the points~\eqref{Eq:D-points}.
In these coordinates, the Bloch variety is defined by $\kappa^2=r^2+s^2$,
showing that it is a double cone, called a \demph{Dirac point}.
These may be observed on the right in Figures~\ref{FigTwo} and~\ref{F:Hex_LandscapeII}.

When $\Delta=0$, the functions $w$ and $z$ are not local coordinates at the critical points, as their
gradients are dependent.
Thus the singularity~\eqref{Eq:Hexagon_ODP} is not an ordinary double point, but has higher order.
Specifically, it has order three when $u\neq v$ and six when $u=v$.
Figure~\ref{F:Hex_LandscapeII} shows three Bloch varieties when $u=v$ with their corresponding parameter values.
%%%%%%%%%%%%%%%%%%%%%%%%%%%%%%%%%%%%%%%%%%%%%%%%%%%%%%%%%%%%%%%%%%%%%%%%%%%%%%%%%%%%%%%%%%%%%%%%%%%%
%%%%%%%%%%%%%%%%%%%%%%%%%%%%%%%%%%%%%%%%%%%%%%%%%%%%%%%%%%%%%%%%%%%%%%%%%%%%%%%%%%%%%%%%%%%%%%%%%%%%

\begin{table}[htb]
\centering
\begin{tabular}{||c | c | c | l||} 
 \hline
 $\Delta$ & nondegenerate &  degenerate & Hessian \\ [0.3ex] 
 & critical points & critical points &  \\ [0.3ex] 
 \hline\hline 
   $< 0$ & 8 & 2 (Dirac)    &\ (1,6,1) \rule{0pt}{12pt}\\ [0.6ex]
   $> 0$ & 8 & 0            &\ (2,4,2) \\ [0.6ex]
   $= 0$ & 6 & 1 (mult.\ 6) &\ (1,4,1) \\ [0.2ex]
  \hline
\end{tabular}
\caption{Critical point summary for the hexagonal lattice\rule{0pt}{14pt} when $u = v$.}
\label{table:hexagonal_third}
\end{table}
%%%%%%%%%%%%%%%%%%%%%%%%%%%%%%%%%%%%%%%%%%%%%%%%%%%%%%%%%%%%%%%%%%%%%%%%%%%%%%%%%%%%%%%%%%%%%%%%%%%%
%

%%%%%%%%%%%%%%%%%%%%%%%%%%%%%%%%%%%%%%%%%%%%%%%%%%%%%%%%%%%%%%%%%%%%%%%%%%%%%%%%%%%%%%%%%%%%%%%%%%%%
\section{Another bipartite graph}\label{S:four}

\begin{figure}[htb]
  \centering
   \begin{picture}(122,92)
     \put(0,0){\includegraphics{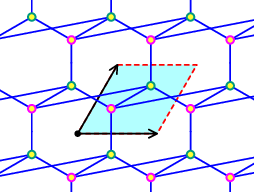}}
   \end{picture}
   \qquad
   \begin{picture}(155,92)(-22,-1)
     \put(-22,0){\includegraphics{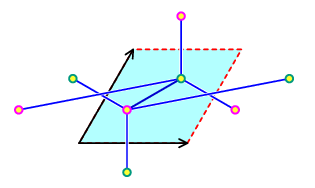}}
     \put(31,31){\small$u$}   \put(95,31){\small$u$}  \put(-22,31){\small$u$}  \put(63,88){\small$u$}
     \put(69,53){\small$v$}   \put( 6,54){\small$v$}  \put(123,54){\small$v$}  \put(37,-2){\small$v$}
     \put(56,14){\small$x$}
     \put(31,61){\small$y$}
     \thicklines 
      \put(63.7,34.7){{\color{White}\vector(-1,1){10}}} \put(64.3,34.3){{\color{White}\vector(-1,1){10}}}
      \put(64.9,34.1){{\color{White}\vector(-1,1){11}}}
     \thinlines
     \put(65.5,29){\small{\color{blue}$a$}} \put(64.4,34.6){\vector(-1,1){10}}
     \put(32,13){\small{\color{blue}$c$}}   \put(68,73){\small{\color{blue}$c$}}
     \put(20,50){\small{\color{blue}$b$}}   \put(81,32){\small{\color{blue}$b$}}
     \put( 9,33){\small{\color{blue}$d$}}   \put(95,51.5){\small{\color{blue}$d$}}
   \end{picture}
   \caption{Bipartite periodic graph with four orbits of edges and a labeling.}
  \label{fig:fourEdgesGraph}
\end{figure}

In this and the next section, we extend some of the analysis of Section~\ref{S:Hexagonal} to two further bipartite graphs
for which the results of Sections~\ref{S:ReciprocalDiscriminants} and~\ref{Sec:bipartite} apply.
Figure~\ref{fig:fourEdgesGraph} shows a bipartite periodic graph $\Gamma$ that has one more orbit of edges than the hexagonal
lattice.
%%%%%%%%%%%%%%%%%%%%%%%%%%%%%%%%%%%%%%%%%%%%%%%%%%%%%%%%%%%%%%%%%%%%%%%%%%%%%%%%%%%%%%%%%%%%%%%%%%%%
%%%%%%%%%%%%%%%%%%%%%%%%%%%%%%%%%%%%%%%%%%%%%%%%%%%%%%%%%%%%%%%%%%%%%%%%%%%%%%%%%%%%%%%%%%%%%%%%%%%%

For this graph and these parameters, we have the Floquet matrix 
\begin{equation}
  \label{Eq:four_FloquetMatrix}
    H_p(x,y) \ = \
     \begin{pmatrix}
          u & a+bx^{-1}+cy^{-1}+dx\\
         a +bx+cy+dx^{-1} & v
     \end{pmatrix}\,,
\end{equation}
({\it cf}.~\eqref{eq:Lgeneral}.)
Its dispersion polynomial $D_p(x,y,\lambda)=\det(\lambda\cdot{\rm id} - H_p(x,y))$ is
 \begin{equation}\label{Eq:Dispersion_four}
      (\lambda-u)(\lambda-v)\ -\ (a+bx^{-1}+cy^{-1}+dx)(a+bx+cy+dx^{-1})\,.
 \end{equation}
We show its support and Newton polytope, $\calN$, which has normalized volume 24.
\[
 \left(\begin{array}{rrrrrrrrrrrrr}
     0&0& 0& 1&-1& 0& 0& 1&-1& 1&-1& 2&-2\\
     0&0& 0& 0& 0& 1&-1&-1& 1& 1&-1& 0& 0\\
     2&1& 0& 0& 0& 0& 0& 0& 0& 0& 0& 0& 0\end{array}\right)
   \qquad
   \raisebox{-25pt}{\begin{picture}(130,60)(-20,0)
       \put(-18,30){$\calN=$}
       \put(0,0){\includegraphics[height=60pt]{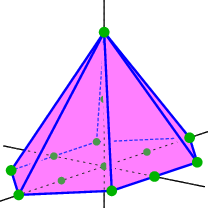}}
       \end{picture}}
\]

Using computer algebra software, we eliminate the variables $x$ and $y$ from the critical point
equations~\eqref{eq:critpoints}.
This gives a degree 14 polynomial in $\lambda$ whose coefficients depend on the parameters $p$.
It factors into two linear, four quadratic, and one quartic polynomial,
\begin{equation} \label{eq:7factors} \begin{matrix}
    \lambda-u\\
    \lambda-v \\
(\lambda-u)(\lambda-v)\ -\ (a+(b+d)+c)^2\\
(\lambda-u)(\lambda-v)\ -\ (a+(b+d)-c)^2\\
(\lambda-u)(\lambda-v)\ -\ (a-(b+d)+c)^2\\
(\lambda-u)(\lambda-v)\ -\ (a-(b+d)-c)^2\\
16b^2d^2(\lambda-u)^2(\lambda-v)^2  \ +\ 8bd(G - 8bc^2d)(\lambda-u)(\lambda-v)\ +\ G^2 \,.
\end{matrix}
\end{equation}
The linear and quadratic factors are those given by Proposition~\ref{prop:obviousfactors}.
Here, \defcolor{$G$} is the quartic
 \begin{equation}~\label{Eq:G-discr}
   G\ \vcentcolon=\  (b-d)^2(a^2-4bd) +4bc^2d\ \ =\ \ 
     a^2(b+d)^2 - 4bd\bigl( a^2+(b-d)^2 -c^2\bigr)\,.
   \vspace{5pt}   %%%%   Funny vertical spacing
 \end{equation}
Restricting the reciprocal $A$-discriminant $\Delta_A^\rec$ of Example~\ref{Ex:medDiscriminant} to the coefficients of
$D_p(x,y,\lambda)$~\eqref{Eq:Dispersion_four} gives 

$b^2 \, c^4 \, d^2 \, (\lambda-u)^2 \, (\lambda-v)^2 \, H$, where \defcolor{$H$} is the seventh factor in~\eqref{eq:7factors}.

%%%%%%%%%%%%%%%%%%%%%%%%%%%%%%%%%%%%%%%%%%%%%%%%%%%%%%%%%%%%%%%%%%%%%%%%%%%%%%%%%%%%%%%%%%%%%%%%%%%%
\begin{theorem}\label{T:anotherBipartite}
  Let $\Gamma$ be the periodic graph in Figure \ref{fig:fourEdgesGraph}.
  The degree $13$ polynomial
 \begin{eqnarray*}
   \BD_\Gamma &=& 
   \bigl( a^2(b+d)^2-16b^2d^2\bigr)\; (u-v) \; G  \; \prod_{\pm,\pm}(a\pm(b+d)\pm c)
 \end{eqnarray*}
 is a Bloch discriminant for $\Gamma$.
\end{theorem}
%%%%%%%%%%%%%%%%%%%%%%%%%%%%%%%%%%%%%%%%%%%%%%%%%%%%%%%%%%%%%%%%%%%%%%%%%%%%%%%%%%%%%%%%%%%%%%%%%%%%

%%%%%%%%%%%%%%%%%%%%%%%%%%%%%%%%%%%%%%%%%%%%%%%%%%%%%%%%%%%%%%%%%%%%%%%%%%%%%%%%%%%%%%%%%%%%%%%%%%%%
\begin{theorem}\label{T:Critical_points_four}
  If $\BD_\Gamma\neq 0$ and $(b-d)(4bd-a^2) \neq 0$, then the complex Bloch variety has $24$ nondegenerate critical points.
  They are distributed among the roots of the factors of \eqref{eq:7factors} as follows.
  \begin{enumerate}
    \item Each critical energy from a quadratic factor corresponds to one critical point.
    \item Each critical energy from a linear factor corresponds to four critical points.
    \item Each critical energy from the seventh factor $H$ corresponds to two critical points.
  \end{enumerate}
\end{theorem}
%%%%%%%%%%%%%%%%%%%%%%%%%%%%%%%%%%%%%%%%%%%%%%%%%%%%%%%%%%%%%%%%%%%%%%%%%%%%%%%%%%%%%%%%%%%%%%%%%%%%

We prove this by analyzing each factor of~\eqref{eq:7factors}. Similar to the analysis in Section~\ref{S:Hexagonal}, the critical energies, 
$\lambda$, are the roots of each factor, from which the $(x,y)$-coordinates of the
corresponding critical point(s) are determined.
The parts (1), (2), and (3) are treated in Sections~\ref{S:Four_corner}, \ref{S:Four_linear}, and ~\ref{S:Four_quartic},
respectively, where we also classify the critical points (in $\TT^2$, in $(\RR^\times)^2$, or in neither).
The multiplicities are partially explained by the following observation:
By the involution $(x,y)\mapsto(x^{-1},y^{-1})$, non-corner critical points will occur in pairs.

%%%%%%%%%%%%%%%%%%%%%%%%%%%%%%%%%%%%%%%%%%%%%%%%%%%%%%%%%%%%%%%%%%%%%%%%%%%%%%%%%%%%%%%%%%%%%%%%%%%%
\subsection{Critical points for quadratic factors.}\label{S:Four_corner}
The analysis of the quadratic factors is similar to the hexagonal lattice, including its independence from  $u-v$.
Their roots are 
 \begin{equation}\label{Eq:lambdas_four}
   \lambda\ =\ \frac{1}{2} \Bigl(u+v\  \pm \sqrt{(u-v)^2 + 4(a \pm(b+d) \pm c)^2}\,\Bigr)\,,
 \end{equation}
and the critical points lie above the corresponding corner points $(x,y)=(\pm1,\pm1)$.
The discriminants are sums of squares, % $(u-v)^2 + 4(a \pm(b+d) \pm c)^2$,
which vanish only when $u-v=(a \pm(b+d) \pm c)=0$, which lies in the $\BD_\Gamma$.
If $u = v$, then~\eqref{Eq:lambdas_four} becomes
\[
\lambda\ =\ u\ \pm\ |a \pm(b+d) \pm c|\,.
\]
This yields eight critical points all lying on the real Bloch variety.
The analysis of these critical points when $\BD_\Gamma\neq 0$ is similar to that given in
Proposition~\ref{P:hexagonal_CPA}, if we replace $b$ there by $b{+}d$ here and $\Delta$ by $\prod (a\pm(b{+}d)\pm c)$.
The main difference is in the analysis of the nature of the critical points. 

%%%%%%%%%%%%%%%%%%%%%%%%%%%%%%%%%%%%%%%%%%%%%%%%%%%%%%%%%%%%%%%%%%%%%%%%%%%%%%%%%%%%%%%%%%%%%%%%%%%%
\subsection{Critical points when $\lambda=u$ or $\lambda=v$.}\label{S:Four_linear}
When $\lambda=u$ or $\lambda=v$,
the Fermi variety is the union of two curves defined by the non-diagonal entries of~\eqref{Eq:four_FloquetMatrix} 
and the critical points are where the curves intersect, as the individual curves are smooth.
These are solutions to the system,
 \begin{equation}\label{Eq:Curve_system}
   a+bx^{-1}+cy^{-1}+dx\ =\ a+bx+cy+dx^{-1}\ =\ 0\,.
 \end{equation}
Solving both for $y$, equating, and simplifying gives the reciprocal quartic 
\begin{equation}
  \label{Eq:Four_x}
0\ =\ bdx^2\ +\ (ab+ad) x\ +\ a^2{+}b^2{+}d^2{-}c^2\ +\ (ab+ad)x^{-1}\ +\ bdx^{-2}\,.
 \end{equation}
If $\defcolor{g(t)} \vcentcolon= 4bdt^2\ +\ 2a(b{+}d)t\ +\ a^2{+}(b{-}d)^2{-}c^2$, 
then the polynomial~\eqref{Eq:Four_x} is $g(\frac{x+x^{-1}}{2})$.
The discriminant of the quadratic $g(t)$ is $4G$, where $G$ is the polynomial~\eqref{Eq:G-discr}.
The polynomial~\eqref{Eq:Four_x} has fewer than four distinct roots if $1$ or $-1$ is a root of $g$, or if $g$ has a
double root, that is, when
\begin{equation} \label{eq:Gg}
  G \,\cdot\, g(1)g(-1)\ =\   G\,\cdot\, \prod_{\pm,\pm}(a\pm(b+d)\pm c)
\end{equation}
vanishes.
This accounts for the last five factors of $\BD_\Gamma$.
The factors $(a\pm(b{+}d)\pm c)$ vanish when the two components of the Fermi variety are tangent at the
corresponding corner point $(x,y)=(\pm1,\pm1)$, as we saw with the hexagonal lattice.
Finally, the factor $G$ vanishes when the components have two points of tangency.

We have the following classification of solutions to~\eqref{Eq:Curve_system}.

%%%%%%%%%%%%%%%%%%%%%%%%%%%%%%%%%%%%%%%%%%%%%%%%%%%%%%%%%%%%%%%%%%%%%%%%%%%%%%%%%%%%%%%%%%%%%%%%%%%%
\begin{theorem}\label{Thm:Classification_Four}
  Suppose that $G\,g(1)g(-1)\neq 0$.
  Then~\eqref{Eq:Curve_system} has four solutions in $(\CC^\times)^2$.
 \begin{enumerate}[label=$(\arabic*)$]
  \item  When $G<0$ no solution lies on $(\RR^\times)^2 \bigcup \TT^2$.
  \item  When $G>0$ all solutions lie on $(\RR^\times)^2 \bigcup \TT^2$.
  \item If $g(1)\cdot g(-1)<0$, then $G>0$ and there are two solutions on $(\RR^\times)^2$ and two on $\TT^2$.
  \item If  $G$, $\frac{1}{bd}g(1)$, $\frac{1}{bd}g(-1)$, and $16b^2d^2{-}a^2(b{+}d)^2$ are positive, then all solutions lie on~$\TT^2$.
  \item In all other cases, all solutions lie on $(\RR^\times)^2$.
 \end{enumerate}
\end{theorem}
%%%%%%%%%%%%%%%%%%%%%%%%%%%%%%%%%%%%%%%%%%%%%%%%%%%%%%%%%%%%%%%%%%%%%%%%%%%%%%%%%%%%%%%%%%%%%%%%%%%%

\begin{figure}[htb]

  \centering
   \begin{picture}(240,180)(-5,0)

    \put(0,0){\includegraphics[height=180pt]{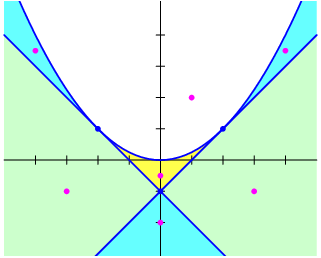}}

    \put( 39,136){\small$B^2{-}4C$}
    \put( 30,  6){\small$1{-}B{+}C$}
    \put(156,  6){\small$1{+}B{+}C$}

    \put(118, 44){\small$(0,-1)$}
    \put( 32, 82){\small$(-2,1)$}  \put(158, 82){\small$(2,1)$}
    \put(209, 72){$B$} \put(100,167){$C$}

    \put( 6,160){\small II}    \put(130,160){\small I}    \put(207,160){\small III}
    \put(22, 38){\small V}   \put(118, 10){\small IV}      \put(190, 38){\small VI}
    \put(59,44){\small VII}  \put( 78,48){\vector(2,1){28}}

   \end{picture}
   %%%%%%%%%%%%%%%%%%%%%%%%%%%%%%%%%%%%%%%%%%%%%%%%%%%%%%%%%%%%%%%%%%%%%%%%%%%%%%%%%%%%%%%%%%%%%%%%%%%%
   \ 
   %%%%%%%%%%%%%%%%%%%%%%%%%%%%%%%%%%%%%%%%%%%%%%%%%%%%%%%%%%%%%%%%%%%%%%%%%%%%%%%%%%%%%%%%%%%%%%%%%%%%
   \raisebox{90pt}{\begin{tabular}{|c|c|c|c|}\hline
      & Point          %& Solutions
       & Summary\\\hline\hline
   I   & $(1, 2)$         %&$-.9052\pm 2.955\sqrt{-1}, -.09478\pm0.3094\sqrt{-1}$\rule{0pt}{12pt}
       & None on $\RR^\times \cup \TT$  \raisebox{-5pt}{\rule{0pt}{17pt}}  \\\hline 
   II  &$(-4, \frac{7}{2})$       %&$.1915,\ .4734,\ 2.112,\ 5.223$
       &  4 on $\RR_>$ \raisebox{-6pt}{\rule{0pt}{18pt}}  \\\hline
   III &$(4, \frac{7}{2})$        %&$-5.223,\ -2.112,\ -.4734,\ -.1915$
       &  4 on $\RR_<$ \raisebox{-6pt}{\rule{0pt}{18pt}}  \\\hline
   IV  &$(0, -2)$         %&$\pm 2.414,\ \pm.4142$
       & 2 on $\RR_<$,  2 on $\RR_>$   \raisebox{-5pt}{\rule{0pt}{17pt}}   \\\hline
    V  &$(-3, -1)$%       &$-.3028\pm.9531\sqrt{-1},\ .1550,\ 6.451$\rule{0pt}{12pt}
       & 2 on $\RR_>$,  2 on $\TT$  \raisebox{-5pt}{\rule{0pt}{17pt}}    \\\hline
   VI  &$(3, -1)$         %&$-6.451, -.1550, .3028\pm .9531\sqrt{-1}$\rule{0pt}{12pt}
       & 2 on $\RR_<$,  2 on $\TT$   \raisebox{-5pt}{\rule{0pt}{17pt}}   \\\hline
   VII &$(0, -\frac{1}{2})$       %&$ \pm .7071 \pm.7071\sqrt{-1}$\rule{0pt}{12pt}
       &    4 on $\TT$ \raisebox{-6pt}{\rule{0pt}{18pt}} \\\hline
\end{tabular}}

  \caption{$\RR^2\smallsetminus\{(B,C)\mid 0=(1+B+C )(1-B+C)(B^2-4C)\}$ has seven connected components.
    We solved~\eqref{Eq:Four_x} at one point in each component.}
   \label{Fig:classification}

\end{figure}

\begin{proof}
  As discussed, the condition $G\,g(1)g(-1)\neq 0$ implies that~\eqref{Eq:Four_x} has four distinct solutions,
  and each gives one solution to~\eqref{Eq:Curve_system}.

  A solution $(x,y)$ to~\eqref{Eq:Curve_system} lies on $\TT^2$ if and only if $x\in\TT$.
  Indeed, as 
\[
   y\ =\ -\frac{a+bx+dx^{-1}}{c} \qquad\mbox{and}\qquad
   y^{-1}\ =\ -\frac{a+bx^{-1}+dx}{c}\ ,
\]
 if $x\in\TT$, then $\overline{x}=x^{-1}$, which implies that $\overline{y}=y^{-1}$, and thus $y\in\TT$.
 The same holds, {\it mutatis mutandis}, if $x\in\RR^\times$.

 To classify the solutions to~\eqref{Eq:Four_x}, let us write
 \begin{equation}\label{Eq:BandC}
  \defcolor{B}\ \vcentcolon=\ \frac{a(b+d)}{2bd} \qquad\mbox{and}\qquad
  \defcolor{C}\ \vcentcolon =\ \frac{a^2 + (b-d)^2 - c^2}{4bd}\,.
 \end{equation}
 Then $g(t) = 4bd(t^2+Bt+C)$ and $G=4b^2d^2(B^2-4C)$, so that 
\[
   g(1) g(-1)\, G\ =\ 64b^4d^4\,( 1+B+C )\,(1-B+C)\,(B^2-4C)\,.
\]
 Let $\defcolor{\mu}\vcentcolon=\frac{1}{2}(-B \pm (B^2-4C)^{1/2})$ be a root of $g$.
Then the solutions $x$ to~\eqref{Eq:Four_x} are
 \begin{equation}\label{Eq:x-solutionsFrom-mu}
   x\ =\ \mu \pm (\mu^2 -1 )^{1/2}
   \qquad\mbox{ with }\qquad x+ x^{-1}\ =\ 2\mu\,.
 \end{equation}
Thus $x\in \RR^\times\cup \TT$ if and only if $\mu\in\RR$, which explains the first two claims.

For the other claims, we solve~\eqref{Eq:Four_x} for one point in
each connected component of $\RR^2$ where $0\neq(1+B+C)(1-B+C)(B^2-4C)$.
This is summarized in Figure~\ref{Fig:classification}.  
We have  $16b^2d^2-a^2(b+d)^2>0$ when $|B|<2$, which separates region VII from regions II and III.
\end{proof}
%%%%%%%%%%%%%%%%%%%%%%%%%%%%%%%%%%%%%%%%%%%%%%%%%%%%%%%%%%%%%%%%%%%%%%%%%%%%%%%%%%%%%%%%%%%%%%%%%%%%

%%%%%%%%%%%%%%%%%%%%%%%%%%%%%%%%%%%%%%%%%%%%%%%%%%%%%%%%%%%%%%%%%%%%%%%%%%%%%%%%%%%%%%%%%%%%%%%%%%%%
%%%%%%%%%%%%%%%%%%%%%%%%%%%%%%%%%%%%%%%%%%%%%%%%%%%%%%%%%%%%%%%%%%%%%%%%%%%%%%%%%%%%%%%%%%%%%%%%%%%%

%%%%%%%%%%%%%%%%%%%%%%%%%%%%%%%%%%%%%%%%%%%%%%%%%%%%%%%%%%%%%%%%%%%%%%%%%%%%%%%%%%%%%%%%%%%%%%%%%%%%
\subsection{Critical points for the seventh factor.}\label{S:Four_quartic}
The seventh factor $H$ in~\eqref{eq:7factors} is a quadratic in
$\defcolor{\nu}\vcentcolon= 4bd(\lambda-u)(\lambda-v)$, namely
\begin{equation}\label{Eq:nu-equation}
  H\ =\ \nu^2\ +\ 2(G-8bc^2d) \nu\ +\ G^2\,.
\end{equation}
This quadratic has discriminant $256b^2c^2d^2(b-d)^2(1 -\frac{a^2}{4bd})$ and roots
\begin{eqnarray}
  \nu &=& 8bc^2d-G\ \pm\ 8bcd(b-d)\bigl( 1-\tfrac{a^2}{4bd}\bigr)^{1/2}\,.\nonumber\\
      &=& 4bd\left( c\ \pm\ (b-d)\bigl(1-\tfrac{a^2}{4bd}\bigr)^{1/2}\right)^2\label{Eq:nu-roots}
\end{eqnarray}
These are real exactly when $1-\frac{a^2}{4bd}\geq 0$ and distinct when $(b-d)(1-\frac{a^2}{4bd})$ is non-zero.
Using~\eqref{Eq:nu-roots} to solve $\nu=4bd(\lambda-u)(\lambda-v)$ for $\lambda$ gives,
 \begin{equation}\label{Eq:lambda_quartic}
   \lambda\ =\ \frac{u+v}{2}\ \ \pm\ \ 
   \sqrt{ \Bigl(\frac{u-v}{2}\Bigr)^2\ +\
   \left( c\ \pm\ (b-d)\bigl(1-\tfrac{a^2}{4bd}\bigr)^{1/2}\right)^2}\ .
 \end{equation}
 Note that the discriminant is a sum of squares.

%%%%%%%%%%%%%%%%%%%%%%%%%%%%%%%%%%%%%%%%%%%%%%%%%%%%%%%%%%%%%%%%%%%%%%%%%%%%%%%%%%%%%%%%%%%%%%%%%%%%
\begin{lemma}\label{L:fourTwo}
  When $u\neq v$, the seventh factor $H$ in~\eqref{eq:7factors} has four roots if and only if $(b-d)(1-\frac{a^2}{4bd})$ is non-zero.
  These are real if\/ $1>a^2/4bd$ and complex otherwise.
  When $(b-d)(1-\frac{a^2}{4bd})$ vanishes, it has two real solutions, each of multiplicity two.
\end{lemma}
%%%%%%%%%%%%%%%%%%%%%%%%%%%%%%%%%%%%%%%%%%%%%%%%%%%%%%%%%%%%%%%%%%%%%%%%%%%%%%%%%%%%%%%%%%%%%%%%%%%%

We expect that each critical energy~\eqref{Eq:lambda_quartic} gives rise to two critical points.
Indeed, using interpolation and a Gr\"obner basis calculation, we find the following three polynomials,
%
% This is done in the maple file FranksComputations/Four/Seventh_Factor/Critical_Points_simple.maple 
%
 \begin{equation}\label{Eq:Four_twoEqs}
   \begin{array}{c}
     2bd\,x \ +\ a(b+d)\ +\ 2bd\,x^{-1}\\
     \pm 2b\bigl(1-\tfrac{a^2}{4bd}\bigr)^{1/2}\cdot y\ -\  (2bx+a)\,\rule{0pt}{13pt}\\
     4bd\bigl(1-\tfrac{a^2}{4bd}\bigr)^{1/2} (y+y^{-1})\ \pm\ 2a(b-d)\,.\rule{0pt}{14pt}
   \end{array}
 \end{equation}
Here, the signs $\pm$ are the signs in the intermediate solution $\nu$~\eqref{Eq:nu-roots}.
For a given value of $\nu$, the first gives two values for $x$, each of which gives one value of $y$, by the second equation.
The discriminant of the first is $a^2(b+d)^2-16b^2d^2$, which we saw in Theorem~\ref{Thm:Classification_Four}(4).
This accounts for the first factor in the Bloch discriminant, $\BD_\Gamma$.
When it is positive, the roots lie in $\RR^\times$ and when negative (the condition in Theorem~\ref{Thm:Classification_Four}(4)),
both are in $\TT$.

%%%%%%%%%%%%%%%%%%%%%%%%%%%%%%%%%%%%%%%%%%%%%%%%%%%%%%%%%%%%%%%%%%%%%%%%%%%%%%%%%%%%%%%%%%%%%%%%%%%%
\begin{theorem}\label{Th:5.5}
  If $16b^2d^2<a^2(b+d)^2$, then there are no critical points on the real Bloch variety arising from the seventh factor.

  Suppose that $16b^2d^2>a^2(b+d)^2$.
 \begin{enumerate}[label=$(\arabic*)$]
 \item
     If $b\neq d$, then the seventh factor $H$ of~\eqref{eq:7factors} has four real solutions, giving four distinct real critical
      energies, and each corresponds to two critical points on the real Bloch variety.
  \item
   If $b=d$, then $H$ has two real solutions, each corresponding to four  critical points on the real Bloch variety. 
 \end{enumerate}
\end{theorem}
%%%%%%%%%%%%%%%%%%%%%%%%%%%%%%%%%%%%%%%%%%%%%%%%%%%%%%%%%%%%%%%%%%%%%%%%%%%%%%%%%%%%%%%%%%%%%%%%%%%%
\begin{proof}
  When $16b^2d^2<a^2(b+d)^2$, the coordinate $x$ of a critical point does not lie on $\TT$, which proves
  the first assertion.
  
  Suppose that $16b^2d^2>a^2(b+d)^2$, so that the coordinate $x$ of a critical point lies on $\TT$.
  The expressions~\eqref{Eq:nu-roots} for $\nu$ and~\eqref{Eq:lambda_quartic} require that $\frac{a^2}{4bd}<1$
  for there to be real critical energies.
  Observe that this always holds if $bd<0$.
  If $bd>0$, then $16b^2d^2>a^2(b+d)^2$ implies that
  \[
  \frac{a^2}{4bd}\ <\ \frac{4bd}{(b+d)^2}\ <\ 1\,.
  \]
  The last inequality is the arithmetic-geometric mean inequality.
  
  For (1), suppose that $b\neq d$.
  As $1>\frac{a^2}{4bd}$, the coefficients of the third equation in \eqref{Eq:Four_twoEqs} are real, so that $y\in\TT\cup\RR^\times$.
  As $x\in\TT$, the second equation shows that $y\not\in \RR$, and thus $y\in\TT$.

  For (2), when $b=d$, we may solve these equations to give four critical points
  \[
     \left(-\frac{a}{2b}   \pm \bigl(\tfrac{a^2}{4b^2}-1\bigr)^{1/2}\ ,\ \pm \sqrt{-1} \right)\,.
  \]
  By Lemma~\ref{L:fourTwo} there are two values of $\lambda$, which completes the proof.
\end{proof}
%%%%%%%%%%%%%%%%%%%%%%%%%%%%%%%%%%%%%%%%%%%%%%%%%%%%%%%%%%%%%%%%%%%%%%%%%%%%%%%%%%%%%%%%%%%%%%%%%%%%

A further analysis, which we do not include, shows that the number and type of \textit{complex} critical points changes when $(b-d)(1-\frac{a^2}{4bd})$ vanishes.
For example, when $b=d$, there are critical points ``at $\infty$'', as studied in~\cite{FLSS,FRS25,FS24}. These factors do
not appear in the Bloch discriminant, as they do not affect the real critical points.
They do, however, appear in the statement of Theorem \ref{Thm:Classification_Four}, which is about the \textit{complex}
critical points of the Bloch variety. 

%%%%%%%%%%%%%%%%%%%%%%%%%%%%%%%%%%%%%%%%%%%%%%%%%%%%%%%%%%%%%%%%%%%%%%%%%%%%%%%%%%%%%%%%%%%%%%%%%%%%
\begin{Remark}
  Section~\ref{S:Four_corner} and Theorems~\ref{Thm:Classification_Four} and~\ref{Th:5.5} give a complete description of the number
  and distribution of the critical points on (real) Bloch varieties when the Bloch discriminant $\BD_\Gamma$ is nonzero,
  which we summarize.
  \begin{enumerate}
    \item There are always two critical points above each of the four corner points.
    
    \item The seventh factor $H$ of~\eqref{eq:7factors} gives either eight or zero critical points, depending upon the sign of
          $16b^2d^2-a^2(b{+}d)^2 =  \left(4bd-a(b{+}d)\right)\left(4bd+a(b{+}d)\right)$.

    \item There are either eight, four, or no critical points when $\lambda$ equals one of the potentials $u$ or $v$.

  \end{enumerate}
  There is always a band gap of at least $|u-v|$ and the gap is exactly $|u-v|$ when there are critical points in (3).
\end{Remark}
%%%%%%%%%%%%%%%%%%%%%%%%%%%%%%%%%%%%%%%%%%%%%%%%%%%%%%%%%%%%%%%%%%%%%%%%%%%%%%%%%%%%%%%%%%%%%%%%%%%%

This classification is succinctly expressed in terms of the coordinates $B$ and $C$~\eqref{Eq:BandC}.
The factor $(u-v)$ of the Bloch discriminant is not evident in the $B,C$-plane.
The other factors map to the four lines and a quadric, $1+B+C$, $1-B+C$, $B-2$, $B+2$ and $B^2-4C$, which divide the plane into
13 discriminant chambers.
We show Bloch varieties with positive parameters $a,b,c,d$ in eight of these, as replacing $B$ by $-B$ does not affect 
our classification, but note that $a,b,c,d>0$ implies that $B>0$.

\begin{figure}[htb]

\centering
  \begin{picture}(124,222)
    \put(0,0){\includegraphics{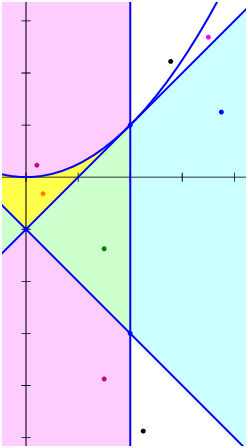}}
    \put(18,205){$C$} \put(108,135){$B$}
    \put(40,190){{\small I}} \put(2,119){{\small II}} \put(40,110){{\small III}} \put(25,10){{\small IV}}
%    \put(22,137){{\small $a$}} 
    \put(80, 200){{\small V}} \put(107, 214){{\small VI}}
    \put(83, 144){{\small VII}} \put(78, 4){{\small VIII}}
  \end{picture}
\quad
\begin{picture}(275,222)(-10,0)
  \put(  0,111){\includegraphics[height=110pt]{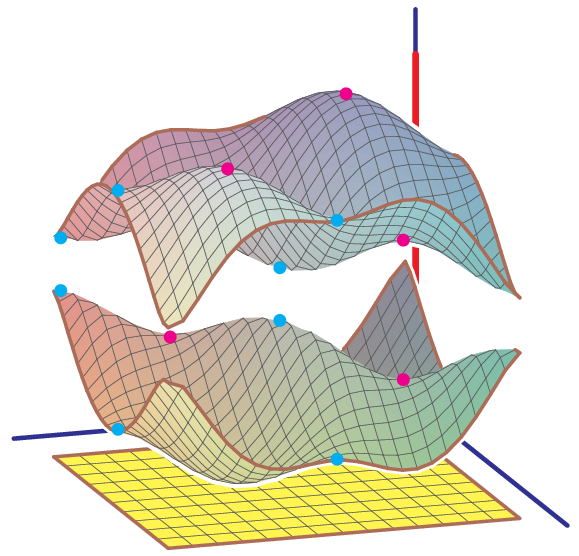}}  % $+-++$ 
  \put(-10,170){{\small I}}
  \put(140,111){\includegraphics[height=110pt]{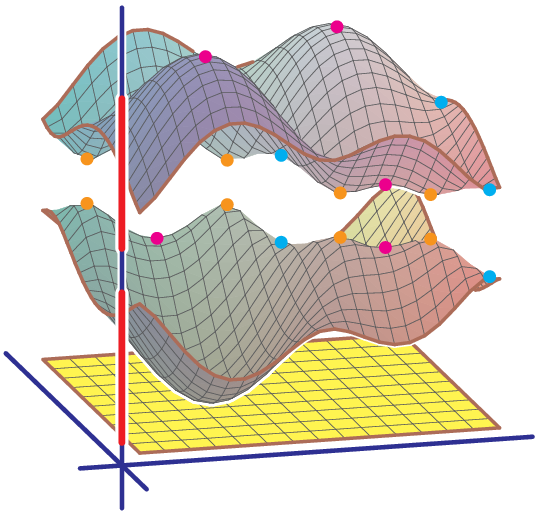}}   % $++++$
  \put(130,170){{\small II}}
  \put(  0,   0){\includegraphics[height=110pt]{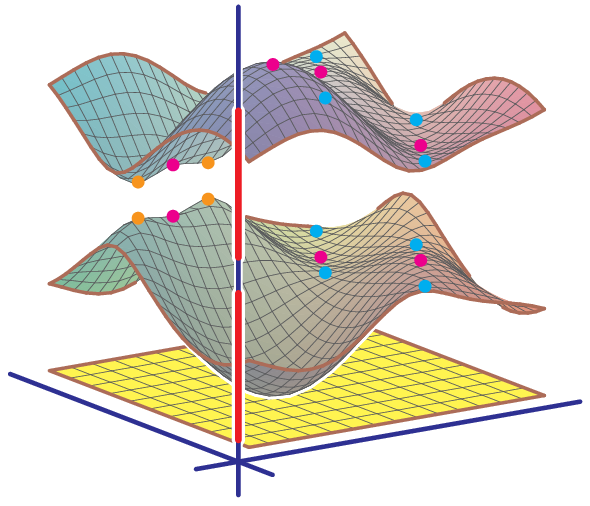}}  % $+++-$
  \put(-10,  60){{\small III}}
  \put(140,   0){\includegraphics[height=110pt]{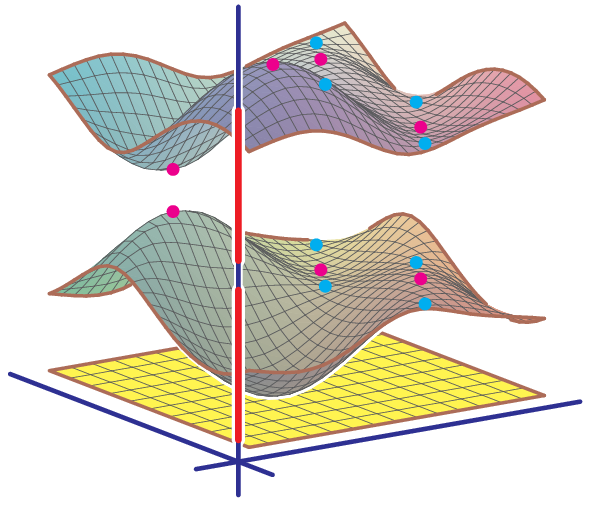}}  % $++--$  
  \put(130,  60){{\small IV}}
\end{picture}
 \caption{Bloch varieties from Chambers I--IV, with 16, 20, or 24 critical points.}
 \label{F:nineI} 
\end{figure}

\begin{figure}[htb]

\centering
\begin{picture}(415,90)
  \put(  0,0){\includegraphics[height=90pt]{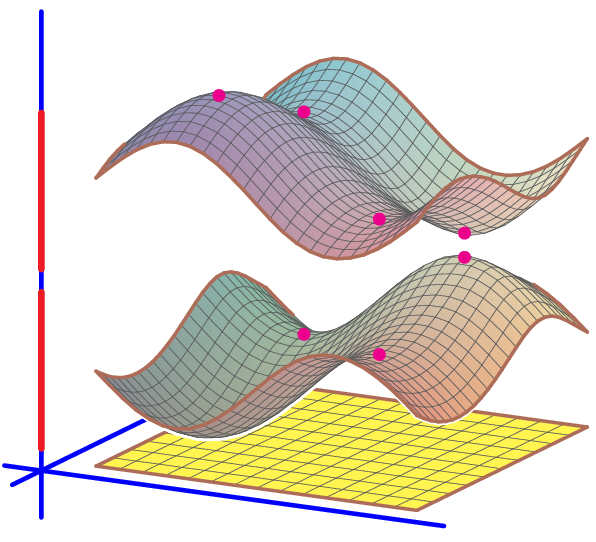}}  % $-+++$ 
     \put( 12,75){{\small V}} 
  \put(105,0){\includegraphics[height=90pt]{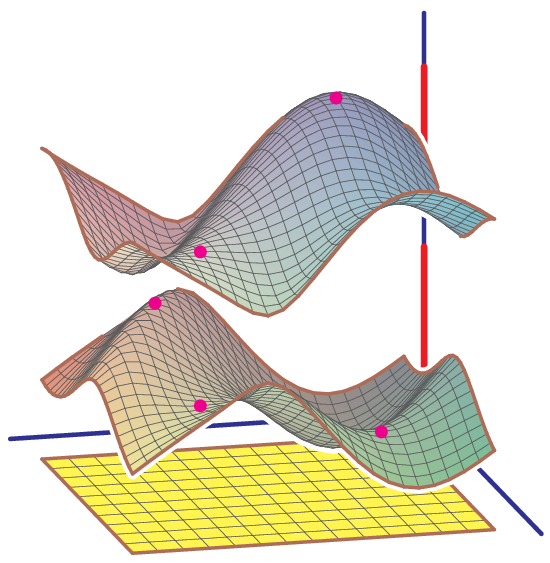}}   % $--++$
     \put(125,75){{\small VI}}
  \put(195,0){\includegraphics[height=90pt]{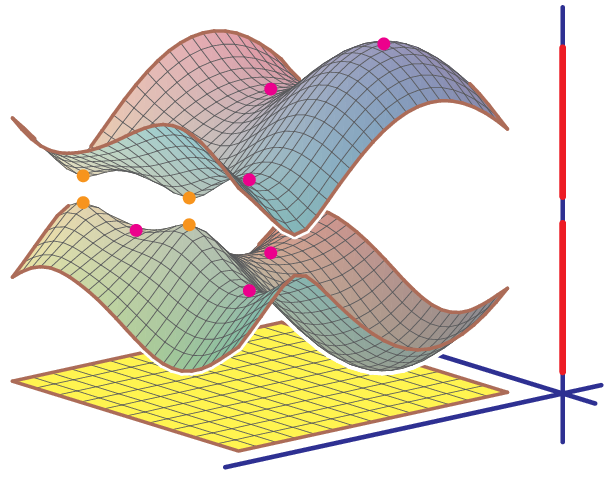}}   % $+-+-$ 
    \put(200,75){{\small VII}} 
  \put(315,0){\includegraphics[height=90pt]{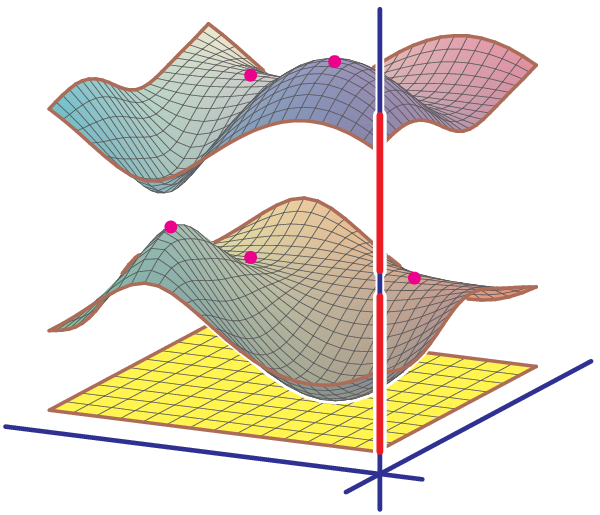}}   %  $-+--$
    \put(315,50){{\small VIII}} 
 \end{picture}

 \caption{Bloch varieties from Chambers V--VIII, with 8 or 12 critical points.}
 \label{F:nineII} 
\end{figure}

Figures~\ref{F:nineI} and~\ref{F:nineII} display eight Bloch varieties, representing one from each chamber.
%%%%%%%%%%%%%%%%%%%%%%%%%%%%%%%%%%%%%%%%%%%%%%%%%%%%%%%%%%%%%%%%%%%%%%%%%%%%%%%%%%%%%%%%%%%%%%%%%%%%
%%%%%%%%%%%%%%%%%%%%%%%%%%%%%%%%%%%%%%%%%%%%%%%%%%%%%%%%%%%%%%%%%%%%%%%%%%%%%%%%%%%%%%%%%%%%%%%%%%%%

%%%%%%%%%%%%%%%%%%%%%%%%%%%%%%%%%%%%%%%%%%%%%%%%%%%%%%%%%%%%%%%%%%%%%%%%%%%%%%%%%%%%%%%%%%%%%%%%%%%%
%%%%%%%%%%%%%%%%%%%%%%%%%%%%%%%%%%%%%%%%%%%%%%%%%%%%%%%%%%%%%%%%%%%%%%%%%%%%%%%%%%%%%%%%%%%%%%%%%%%%

\begin{figure}[htb]
  \centering
  \begin{picture}(160,100)
    \put(  0,  0){\includegraphics[height=110pt]{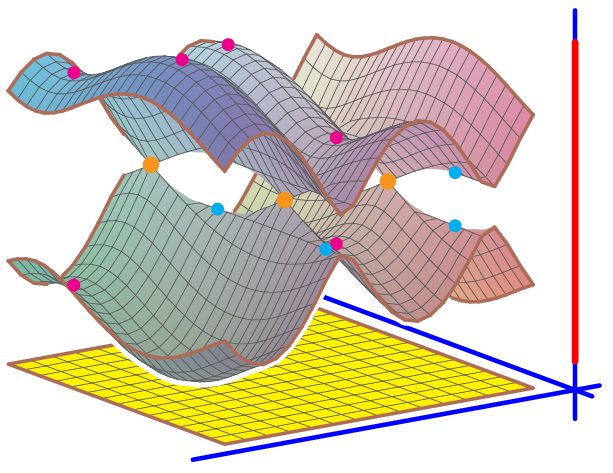}}
    \put( 0,60){{\small II}} 
\end{picture}
\begin{picture}(130,115)
    \put(  0,  0){\includegraphics[height=115pt]{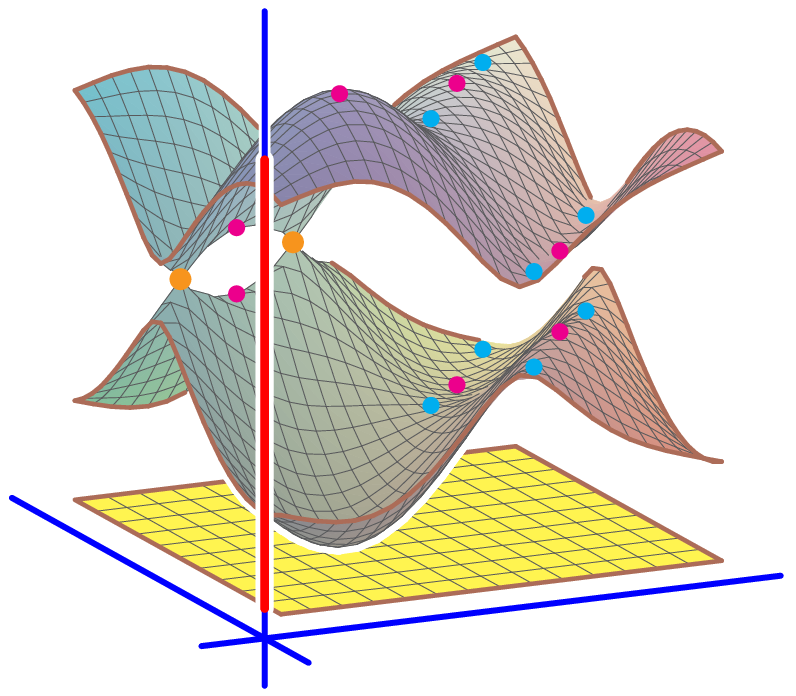}}
    \put( 0,65){{\small III}} 
\end{picture}
\begin{picture}(130,115)
    \put(10,  0){\includegraphics[height=115pt]{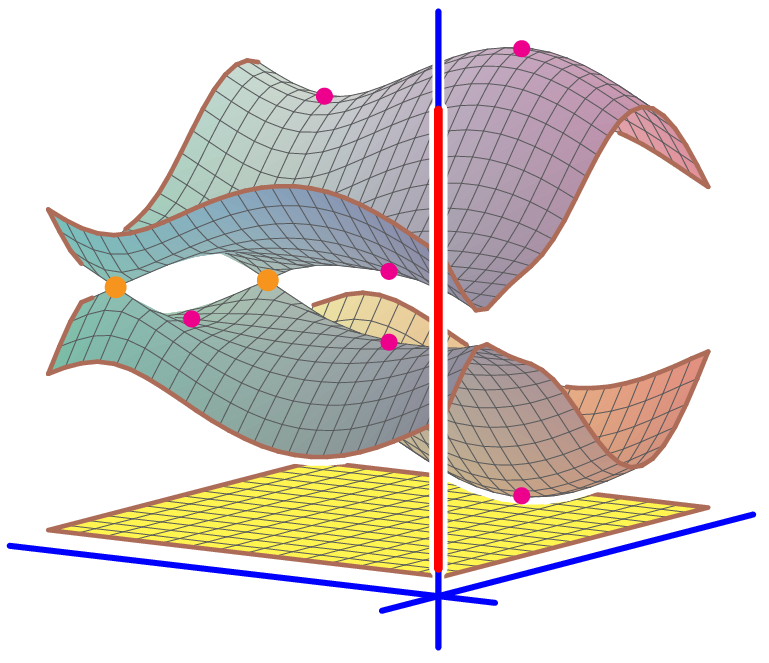}}
    \put( 4, 85){{\small VII}} 
\end{picture}
\caption{Representative Bloch varieties with Dirac points}
\label{F:Four_Dirac}
\end{figure}

Let $G, g$ be as in \eqref{eq:Gg}. When $u=v$ but $(b-d)(1-\frac{a^2}{4bd})\, G \, g(1)g(-1) \neq 0$, the complex Bloch variety has four
ordinary double points.
These give Dirac points only when the parameters $B,C$ lie in one of the Chambers II, III, and VII---in all other
cases, the type and distribution of critical points is qualitatively equivalent to that when $u\neq v$ as in
Figures~\ref{F:nineI} and~\ref{F:nineII}.
There are four Dirac points when $(B,C)$ lies in Chamber II and two for parameters in Chambers III and VII.
Representative Bloch varieties for each are depicted in Figure~\ref{F:Four_Dirac}.
This accounts for the factor $(u-v)$ in the Bloch discriminant.

%%%%%%%%%%%%%%%%%%%%%%%%%%%%%%%%%%%%%%%%%%%%%%%%%%%%%%%%%%%%%%%%%%%%%%%%%%%%%%%%%%%%%%%%%%%%%%%%%%%%
%%%%%%%%%%%%%%%%%%%%%%%%%%%%%%%%%%%%%%%%%%%%%%%%%%%%%%%%%%%%%%%%%%%%%%%%%%%%%%%%%%%%%%%%%%%%%%%%%%%%

%
%%%%%%%%%%%%%%%%%%%%%%%%%%%%%%%%%%%%%%%%%%%%%%%%%%%%%%%%%%%%%%%%%%%%%%%%%%%%%%%%%%%%%%%%%%%%%%%%%%%%
\section{A third bipartite graph}\label{S:five}

\begin{figure}[htb]
  \centering
   \begin{picture}(122,105)(0,-11)
     \put(0,0){\includegraphics{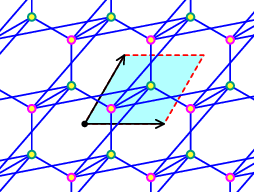}}
     \put(77,54){$W$}
   \end{picture}
   \qquad
   \begin{picture}(153,115)(-22,-10)
     \put(-22,-15){\includegraphics{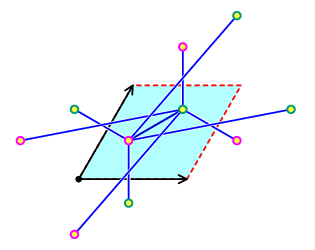}}
     
    \thicklines 
      \put(66.7,35.7){{\color{White}\vector(-1,1){10}}} \put(67.3,36.3){{\color{White}\vector(-1,1){10}}}
      \put(67.9,35.1){{\color{White}\vector(-1,1){11}}}

     \thinlines
     \put(31,31){\small$u$}   \put(110,29){\small$u$}  \put(-21,32){\small$u$}  \put(63,100){\small$u$}
     \put(108,32){\vector(-3,1){13}}     \put(66,98){\vector(0,-1){12}}
     \put(69,54){\small$v$}   \put( 6,54){\small$v$}  \put(121,54){\small$v$}  \put(37,-1){\small$v$}
     \put(58,23){\small$x$}
     \put(33,63){\small$y$}

     \put(68,31){\small{\color{blue}$a$}}   \put(67.6,35.4){\vector(-1,1){10}}
     \put(42,13){\small{\color{blue}$c$}}   \put(59,75 ){\small{\color{blue}$c$}}
     \put(23,11){\small{\color{blue}$e$}}   \put(78,76  ){\small{\color{blue}$e$}}
     \put(20,50){\small{\color{blue}$b$}}   \put(83,32 ){\small{\color{blue}$b$}}
     \put(14,34){\small{\color{blue}$d$}}   \put(93,50.5){\small{\color{blue}$d$}}
   \end{picture}
   \caption{Bipartite periodic graph with five orbits of edges and a labeling.}
  \label{fig:fiveEdgesGraph}
\end{figure}

We partially analyze the critical points of operators for a third  bipartite graph.
Let $\Gamma$ be the bipartite graph in Figure~\ref{fig:fiveEdgesGraph}
%%%%%%%%%%%%%%%%%%%%%%%%%%%%%%%%%%%%%%%%%%%%%%%%%%%%%%%%%%%%%%%%%%%%%%%%%%%%%%%%%%%%%%%%%%%%%%%%%%%%
%%%%%%%%%%%%%%%%%%%%%%%%%%%%%%%%%%%%%%%%%%%%%%%%%%%%%%%%%%%%%%%%%%%%%%%%%%%%%%%%%%%%%%%%%%%%%%%%%%%%
with two vertices in its fundamental domain and five orbits of edges.
For this graph and these parameters, we have the Floquet matrix,
\begin{equation}
  \label{Eq:five_FloquetMatrix}
    H_p(x,y) \ = \
     \begin{pmatrix}
          u & a+bx^{-1}+cy^{-1}+dx+ey\\
         a +bx+cy+dx^{-1}+ey^{-1} & v
     \end{pmatrix}\,,
\end{equation}
with dispersion polynomial $D_p(x,y,\lambda)=\det(\lambda\cdot{\rm id} - H_p(x,y))$ 
 \begin{equation}\label{Eq:Dispersion_five}
      (\lambda-u)(\lambda-v)\ -\ (a+bx^{-1}+cy^{-1}+dx+ey)(a+bx+cy+dx^{-1}+ey^{-1})\,.
 \end{equation}
We show its support and Newton polytope, $\calN$, which has normalized volume 32.
{\small\begin{equation}\label{Eq:Third_support}
 \left(\begin{array}{rrrrrrrrrrrrrrr}
     0\!&0\!& 0\!& 1\!\!&-1\!& 0\!\!& 0\!& 1\!&-1\!& 1\!\!&-1\!& 2\!\!&-2\!&0\!\!& 0\\
     0\!&0\!& 0\!& 0\!\!& 0\!& 1\!\!&-1\!&-1\!& 1\!& 1\!\!&-1\!& 0\!\!& 0\!&2\!\!&-2\\
     2\!&1\!& 0\!& 0\!\!& 0\!& 0\!\!& 0\!& 0\!& 0\!& 0\!\!& 0\!& 0\!\!& 0\!&0\!\!& 0\end{array}\right)
   \quad
   \raisebox{-22pt}{\begin{picture}(120,55)(-40,0)
       \put(-25,30){$\calN=$}
       \put(-10,0){\includegraphics[height=55pt]{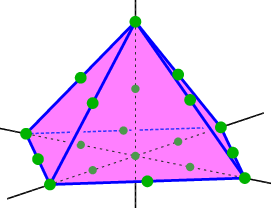}}
       \end{picture}}
\end{equation}}

Unlike the previous two graphs, we are unable to compute any  Gr\"obner basis for its critical point
equations~\eqref{eq:critpoints},
where we treat the parameters $a,b,c,d,e,u,v$ as indeterminates.  This is
a prerequisite for symbolic elimination to obtain a polynomial in $\lambda$ and the parameters,
which we call the $\lambda$-eliminant.
Its roots are the critical energies, the values of $\lambda$ whose corresponding Fermi variety is singular, by
Lemma~\ref{L:singularFermi}.

As a polynomial in $(x,y)$, the support $A$ of the dispersion polynomial is the base of the pyramidal Newton polytope
in~\eqref{Eq:Third_support}.
The corresponding discriminant $\Delta_r(A)$ was described in Example~\ref{Ex:bigDiscriminant}.   
Evaluating $\Delta_r(A)$ at the $(x,y)$-coefficients of the dispersion polynomial gives a polynomial with nine factors.
By the discussion in the previous paragraph, the factors involving $\lambda$ divide the $\lambda$-eliminant.
Four
factors
are squares,
\[
    (be-cd)^2\,,\quad (bc-de)^2\,,\ \quad\   (\lambda-u)^2\,, \quad \mbox{ and }\quad (\lambda-v)^2\,,
\]
there are four irreducible quadratics,
 \begin{equation}\label{Eq:five_quadratics}
   (\lambda-u)(\lambda-v)\ -\ \bigl( a \pm (b+d) \pm (c+e) \bigr)^2\,,
 \end{equation}
and one factor of degree sixteen with 2418 terms.
Let us write \defcolor{$H$} for this last factor.

By inspection (and through hindsight from the form of the dispersion polynomial~\eqref{Eq:Dispersion_five}) the $\lambda$-eliminant is a
polynomial in $\defcolor{\Lambda}\vcentcolon=(\lambda-u)(\lambda-v)$.
Let
\[
  \defcolor{F}\ \vcentcolon=\ x^2y^2\bigl( \Lambda\ -\ (a+bx^{-1}+cy^{-1}+dx+ey)(a+bx+cy+dx^{-1}+ey^{-1})\bigr)\,,
\]
be the numerator of the dispersion polynomial~\eqref{Eq:Dispersion_five} after this substitution.
The critical point ideal becomes the ideal $I\subset\QQ[w,z,x,y,\, \Lambda, u,v,a,b,c,d,e]$ with generators
\[
   F\,,\ \frac{\partial}{\partial x} F\,,\ \frac{\partial}{\partial y} F\,,\ xw-1\,,\ yz-1\,,
\]
Both maple and Singular take about an hour to compute a degree reverse lexicographic Gr\"obner basis for $I$ whose penultimate
member happens to lie in $\QQ[\Lambda, u,v,a,b,c,d,e]$.
This is essentially Algorithm~\ref{alg:discr}, and it stops in Step 4.
That polynomial has six factors.
The first five have low degree in $\Lambda$,
 \begin{equation}\label{Eq:low_Lambda}
   \Lambda\,,\  \  \Lambda - \bigl( a \pm (b+d) \pm (c+e) \bigr)^2\,,
 \end{equation}
and the last, \defcolor{$G(\Lambda)$}, is a quartic in $\Lambda$ with 843 terms and degree 16 in the parameters with a highly
structured constant term,
  \begin{multline}\label{Eq:five_Big}
   \qquad 256 b^2c^2d^2e^2\Lambda^4\ +\ 128 bcde\Lambda^3( a^2b^2ce\ +\ \dotsb\,) +\ \Lambda^2(\,\dotsb\,)\ +\ \Lambda(\,\dotsb\,)\\
   +\  \underline{\left( (a^2 - 4(bd + ce) )^2 - 64bcde\right)}\\
   \times \Bigl(a^2((c-e)^2(b-d)^2 ) + 4 ( (b-d)^2-(c-e)^2 ) (bc-ed)(be-cd) \Bigr)^2 \,. \qquad 
 \end{multline}
Below, we write $D_1 \cdot D_2^{2}$ for the constant term $G(0)$, where $D_1$ is the underlined factor in the display above. Substituting $(\lambda-u)(\lambda-v)$ for $\Lambda$ into the factors~\eqref{Eq:low_Lambda} gives the linear~factors
\[
(\lambda-u) \quad\mbox{ and }\quad (\lambda-v)
\]
and the quadratic factors~\eqref{Eq:five_quadratics}.
This substitution in the polynomial $G(\Lambda)$ gives $H(\lambda)$.
The product of the quadratic factors~\eqref{Eq:five_quadratics}, the linear factors $(\lambda-u)(\lambda-v)$, and the octic polynomial
$H(\lambda)$ gives the generic $\lambda$-eliminant of the critical point equations~\eqref{eq:critpoints}.

Even though $H(\lambda)$ is solvable (as $G(\Lambda)$ is a quartic), its complexity precludes a detailed study of the critical points
for this family of Bloch varieties. {However, our methods do provide a Bloch discriminant.}

\begin{theorem} \label{Thm:Bloch_five}
Let $\Gamma$ be the periodic graph in Figure \ref{fig:fiveEdgesGraph}. The polynomial 
  \[
    {\rm BD}_\Gamma \  = \  (u-v) \cdot {\rm disc}_\lambda H(\lambda) \cdot D_1 \cdot D_2
    \prod_{\pm,\pm}  D_{\pm , \pm}  \, \prod_{\pm , \pm} G( D_{\pm , \pm} ^2)\, ,
   \]
where $D_{\pm , \pm} = a \pm (b+d) \pm (c+e)$, is a Bloch discriminant for $\Gamma$.
\end{theorem}

The factor ${\rm disc}_\lambda H(\lambda)$ in ${\rm BD}_\Gamma$ is the discriminant polynomial which vanishes when 
$H(\lambda)$ has a multiple root.
Since $H$ is a monic polynomial of degree eight in $\lambda$, its discriminant can be written as the determinant of
the $15\times 15$ Sylvester matrix for $H(\lambda), H'(\lambda)$,
whose entries are polynomials in $a,b,c,d,e,u,v$.
This gives an effective {way to evaluate} ${\rm disc}_\lambda H(\lambda)$. We will prove Theorem \ref{Thm:Bloch_five} below.

A symbolic computation at numerical parameters $(a,b,c,d,e,u,v)=(2,1,1,2,3,1,0)$ gives one critical point for each of the
two roots of each quadratic factor~\eqref{Eq:five_quadratics}, and these are above the corresponding corner point $(\pm1,\pm1)$.
For each of the linear factors, there are four critical points, and for each of the eight distinct roots of $H$, there are two critical
points.
This number, 32, of critical points, is the volume of the Newton polytope in~\eqref{Eq:Third_support}, which is an upper bound for the
number of critical points.
As explained in~\cite{DKS,FS24}, this computation establishes the spectral edges nondegeneracy conjecture for the graph of
Figure~\ref{fig:fiveEdgesGraph}.
It also implies the following:

%%%%%%%%%%%%%%%%%%%%%%%%%%%%%%%%%%%%%%%%%%%%%%%%%%%%%%%%%%%%%%%%%%%%%%%%%%%%%%%%%%%%%%%%%%%%%%%%%%%%
\begin{theorem}\label{T:Critical_points_five}
  For generic parameters, the complex Bloch variety has $32$ nondegenerate critical points.
  They are distributed among the roots of the factors as follows.
  \begin{enumerate}
    \item Each critical energy from a quadratic factor corresponds to one critical point.
    \item Each critical energy from a linear factor corresponds to four critical points.
    \item Each critical energy from the factor $H$ corresponds to two critical points.
  \end{enumerate}
\end{theorem}
%%%%%%%%%%%%%%%%%%%%%%%%%%%%%%%%%%%%%%%%%%%%%%%%%%%%%%%%%%%%%%%%%%%%%%%%%%%%%%%%%%%%%%%%%%%%%%%%%%%%

%%%%%%%%%%%%%%%%%%%%%%%%%%%%%%%%%%%%%%%%%%%%%%%%%%%%%%%%%%%%%%%%%%%%%%%%%%%%%%%%%%%%%%%%%%%%%%%%%%%%
\subsection{Critical points for quadratic factors.}\label{S:Five_corner}
The analysis of the critical points corresponding to the quadratic factors
\[
   (\lambda-u)(\lambda-v)\ -\ \left( a \pm (b+d) \pm (c+e) \right)^2\,,
\]
which lie above the corresponding corner points $(\pm 1,\pm 1)$ is identical to that in~\eqref{S:Four_corner}, here
replacing $\pm c$ by $\pm(c+e)$.

%%%%%%%%%%%%%%%%%%%%%%%%%%%%%%%%%%%%%%%%%%%%%%%%%%%%%%%%%%%%%%%%%%%%%%%%%%%%%%%%%%%%%%%%%%%%%%%%%%%%
\subsection{Critical points when $\lambda=u$ or $\lambda=v$.}\label{S:Five_linear}
\begin{figure}[htb]
 \centering 
   \includegraphics[height=150pt]{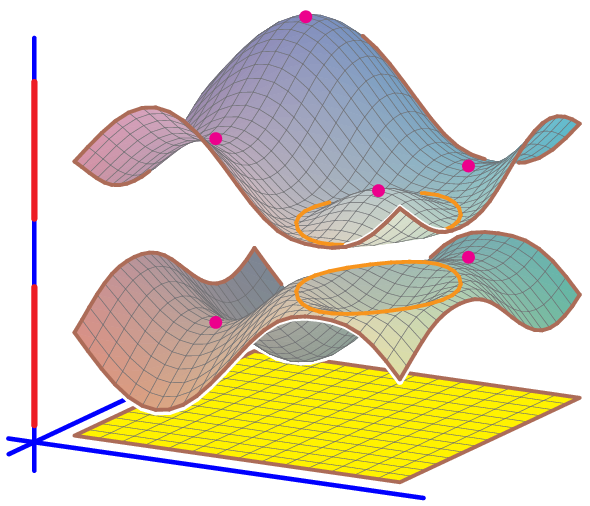}\qquad\quad
   \includegraphics[height=150pt]{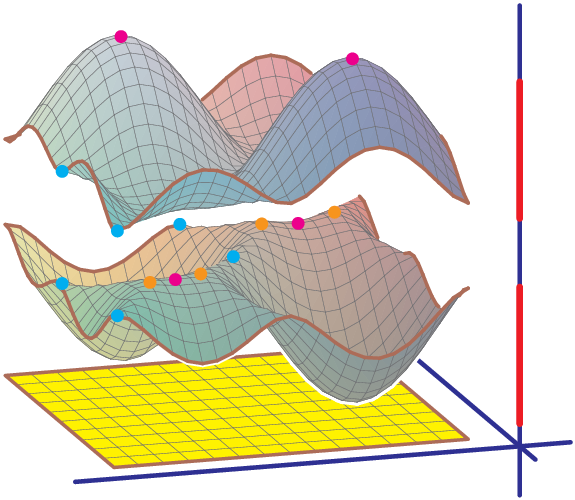}

   \caption{A Bloch variety with a curve of critical points and one with 24 critical points.}
   \label{F:FiveBlochVars}
\end{figure}

The two cases are identical when $u\neq v$.
Suppose that $\lambda=u\neq v$.
Two new phenomena occur for the graph $\Gamma$ which do not occur in~\S~\ref{S:Four_linear} or~\ref{S:Four_quartic}.
If $b=d$ and $c=e$, then the Fermi curve becomes
 \begin{equation}\label{Eq:sing_Fermi}
   \left( a + b(x+ x^{-1}) + c(y+y^{-1}) \right)^2\ =\ 0\,,
 \end{equation}
and there will be a 1-dimensional critical point locus when~\eqref{Eq:sing_Fermi} has solutions in $\TT^2$.
This is essentially the example of Filonov-Kachkovskiy~\cite{FK18}.
We display such a Bloch variety in Figure~\ref{F:FiveBlochVars}, along with its curve of critical points.
It has parameters $(a,b,c,d,e,u,v)=(5,2,3,2,3,6,0)$.
Notice that the polynomial $D_2$ vanishes when $ b = d$ and $c = e$.
This is the first new phenomenon.
 
When $b\neq d$ or $c\neq e$, the singular points of the Fermi curve
\begin{equation} \label{eq:fermicurve_five}
   \{(x,y)\in(\CC^\times)^2\mid ( a + bx+cy + dx^{-1} + ey^{-1}) ( a + bx^{-1}+cy^{-1} + dx + ey)=0\}
\end{equation}
are isolated.
As before, these include the points where both factors vanish.
Both factors have the same Newton polygon
\[
\includegraphics{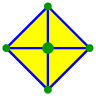}
\]
which has normalized area 4, so for general $a,b,c,d,e$, there will be four such points.

In the context of our paper, we are interested in characterizing parameters which are not \textit{general} in this sense. Using variable elimination, we compute the discriminant for 
\[ a + bx+cy + dx^{-1} + ey^{-1} \ = \ a + bx^{-1}+cy^{-1} + dx + ey \ = \ 0\,, \] 
i.e., a polynomial which vanishes when these equations have a singular solution. The result~is 
\begin{equation} \label{eq:mixeddisc} D_2 \cdot \prod a \pm (b+d) \pm (c+e)\,. \end{equation}
This gives us all the tools we need to prove Theorem \ref{Thm:Bloch_five}. 
\begin{proof}[Proof of Theorem \ref{Thm:Bloch_five}]
  Suppose that the $\lambda$-eliminant has $18$ distinct roots (viewed as a polynomial in $\lambda$).
  We claim that, provided ${\rm BD}_\Gamma \neq 0$, each root $\lambda^*$ of the $\lambda$-eliminant corresponds to the
  expected number of critical points from Theorem \ref{T:Critical_points_five}. This is clearly true if $\lambda^*$ corresponds to only one
  critical point, i.e., $\lambda^*$ is a root of a quadratic factor. For $\lambda^* = u$ or $\lambda^* = v$, the four critical points are
  distinct by the non-vanishing of \eqref{eq:mixeddisc}. Finally, for a root $H(\lambda^*) = 0$, the two corresponding critical points
  collide if and only if they lie at a corner point. But then $\lambda^*$ is also a root of one of the quadratic factors, which is a
  contradiction because {$G((a \pm (b+d) \pm(c+e))^2) \neq 0$} for each choice of signs. 
  
  Now suppose that the $\lambda$-eliminant does not have the expected number of $18$ distinct roots and ${\rm BD}_\Gamma \neq 0$.
  In that case, two of the quadratic factors of the $\lambda$-eliminant must have a common root.
  This is because all other discriminants and pairwise resultants of factors of the $\lambda$-eliminant divide ${\rm BD}_\Gamma$. 
  Such a common root, however, corresponds to the critical energy of distinct critical points whose $(x,y)$ coordinates are
  distinct corner points. Hence, since ${\rm BD}_\Gamma \neq 0$, the Bloch variety has 32 nondegenerate critical points. 
\end{proof}

Each factor in \eqref{eq:fermicurve_five} defines a curve in $(\CC^\times)^2$ of arithmetic genus 1
(given by the number of internal integer points in the Newton polygon~\cite{Baker1894}).
However, when the $A$-discriminant
 \begin{equation}\label{Eq:diamondDiscriminant}
   D_1 \, = \, (a^2-4bd-4ce)^2\ -\ 64bcde
 \end{equation}
of either factor of $ ( a + bx+cy + dx^{-1} + ey^{-1}) ( a + bx^{-1}+cy^{-1} + dx + ey)$ 
vanishes, each factor defines a rational curve with a unique singular point.
This is the second new phenomenon.
This discriminant is a factor of the constant term in the quartic $G(\Lambda)$ (it is the underlined term in~\eqref{Eq:five_Big}).
Thus when it vanishes, the product $(\lambda-u)(\lambda-v)$ of the two linear factors divides $H(\lambda)$.
Let  $K(\lambda)$ be the other factor of $H(\lambda)$, a sextic.
We give a version of Theorem~\ref{T:Critical_points_five}, describing what happens when the
discriminant~\eqref{Eq:diamondDiscriminant} vanishes.

%%%%%%%%%%%%%%%%%%%%%%%%%%%%%%%%%%%%%%%%%%%%%%%%%%%%%%%%%%%%%%%%%%%%%%%%%%%%%%%%%%%%%%%%%%%%%%%%%%%%
\begin{theorem}\label{T:noSingPts}
  Suppose that the discriminant~\eqref{Eq:diamondDiscriminant} vanishes, but that the parameters are otherwise general.
  Then the complex Bloch variety has $32$ nondegenerate critical points.
  They are distributed among the roots of the factors as follows.
  \begin{enumerate}
    \item Each critical energy from a quadratic factor corresponds to one critical point.
    \item Each critical energy from a linear factor corresponds to six critical points, two of which
    arise from the singular points of the irreducible components of the Fermi curve.
    \item Each critical energy from the sextic factor $K$ corresponds to two critical points.
  \end{enumerate}
  When the parameters $a,b,c,d,e$ are real, at most four of the critical points in (2) lie on the real Bloch variety.
\end{theorem}
%%%%%%%%%%%%%%%%%%%%%%%%%%%%%%%%%%%%%%%%%%%%%%%%%%%%%%%%%%%%%%%%%%%%%%%%%%%%%%%%%%%%%%%%%%%%%%%%%%%%
\begin{proof}
  When $(a,b,c,d,e,u,v)=(2,1,1,4,9,1,0)$, the  discriminant~\eqref{Eq:diamondDiscriminant} vanishes and a numerical calculation
  verifies that Statements (1)--(3) hold, which implies the claim about the complex Bloch variety.

  For the last assertion, note that when $a,\dotsc,e$ are real, each curve defined by a factor
  of
  \[
     ( a + bx+cy + dx^{-1} + ey^{-1}) ( a + bx^{-1}+cy^{-1} + dx + ey)
  \] 
  is invariant under ordinary
  complex conjugation, which implies that its singular point is also invariant, and is hence a point in $(\RR^\times)^2$.
  As the parameters are otherwise general, these singular points do not occur at corner points, and hence not on the real Bloch variety.
\end{proof}
%%%%%%%%%%%%%%%%%%%%%%%%%%%%%%%%%%%%%%%%%%%%%%%%%%%%%%%%%%%%%%%%%%%%%%%%%%%%%%%%%%%%%%%%%%%%%%%%%%%%

%%%%%%%%%%%%%%%%%%%%%%%%%%%%%%%%%%%%%%%%%%%%%%%%%%%%%%%%%%%%%%%%%%%%%%%%%%%%%%%%%%%%%%%%%%%%%%%%%%%%
\begin{Remark}\label{R:MaxCP}
  While for general parameters, the complex Bloch variety has 32 critical points, the real Bloch variety
  typically has fewer critical points.
  In both Sections~\ref{S:Hexagonal} and~\ref{S:four}, we exhibited parameters for which the real Bloch variety has the
  maximal number of critical points (12 and 24, respectively).
  In contrast, for the graph here, we have  been unable to find parameters with more than 24 critical points on their
  real Bloch variety.
  We computed the number of real critical points for  $1.5$ million random choices of parameters
  $(a,b,c,d,e,u,v)\in\{1,\dotsc,12\}^{7}$.
  For every even number $k$ between 8 and 24, this search found parameters whose real Bloch variety had $k$ critical
  points.
  For $k\in\{10,14,18,22\}$ these included singular critical points above a corner point and the $\lambda$-eliminant did
  not have 18 distinct roots. 
    For other values of $k$, the $\lambda$-eliminant had 18 distinct roots.
  It found no real Bloch varieties with more than 24 critical points.
  Here, 8 is the lower bound given by the corner critical points.
  On the right in Figure~\ref{F:FiveBlochVars}, we display the real Bloch variety for the parameters
  $(1, 4, 10, 11, 8, 1, 15)$, which has the maximum observed number of 24 critical points.
\end{Remark}
%%%%%%%%%%%%%%%%%%%%%%%%%%%%%%%%%%%%%%%%%%%%%%%%%%%%%%%%%%%%%%%%%%%%%%%%%%%%%%%%%%%%%%%%%%%%%%%%%%%%

%%%%%%%%%%%%%%%%%%%%%%%%%%%%%%%%%%%%%%%%%%%%%%%%%%%%%%%%%%%%%%%%%%%%%%%%%%%%%%%%%%%%%%%%%%%%%%%%%%%%
%%%%%%%%%%%%%%%%%%%%%%%%%%%%%%%%%%%%%%%%%%%%%%%%%%%%%%%%%%%%%%%%%%%%%%%%%%%%%%%%%%%%%%%%%%%%%%%%%%%%

%%%%%%%%%%%%%%%%%%%%%%%%%%%%%%%%%%%%%%%%%%%%%%%%%%%%%%%%%%%%%%%%%%%%%%%%%%%%%%%%%%%%%%%%%%%%%%%%%%%%

\bibliographystyle{abbrv}

\end{document}